\documentclass{amsart}
\usepackage{amsmath,amssymb,amscd,latexsym}
\usepackage[all]{xy}

\newcommand{\C}{{\mathbb{C}}}

\newcommand{\N}{{\mathbb{N}}}

\newcommand{\R}{{\mathbb{R}}}

\newcommand{\Bh}{{\mathcal B}}
\newcommand{\Ch}{{\mathcal C}}

\newcommand{\Fh}{{\mathcal F}}
\newcommand{\Gh}{{\mathcal G}}

\newcommand{\Zh}{{\mathcal Z}}

\newcommand{\Aff}{\mathrm{Aff}}

\newcommand{\be}{\mathbf{1}}

\newcommand{\bj}{\bar{\jmath}}

\newcommand{\dist}{\mathrm{dist}}
\newcommand{\dr}{\mathrm{dr}\,}

\newcommand{\halb}{\frac{1}{2}}

\newcommand{\id}{\mathrm{id}}
\newcommand{\ord}{\mathrm{ord}\,}

\newcounter{number}[section]

\newenvironment{nummer}{\refstepcounter{number}{\noindent\arabic{section}.\arabic{number}}}{}

\newcommand{\bn}{\noindent \begin{nummer} \rm}
\newcommand{\en}{\end{nummer}}

\newenvironment{ntheorem}{\noindent {\sc Theorem:} \it}{}
\newenvironment{nlemma}{\noindent {\sc Lemma:} \it}{}
\newenvironment{nprop}{\noindent {\sc Proposition:} \it}{}
\newenvironment{ndefn}{\noindent {\sc Definition:} \it}{}
\newenvironment{ncor}{\noindent {\sc Corollary:} \it}{}

\newenvironment{nremarks}{\noindent {\sc Remarks:}}{}
\newenvironment{nexamples}{\noindent {\sc Examples:} }{}

\newenvironment{nnotation}{\noindent {\sc Notation:} }{}

\newenvironment{nproof}{\noindent {\sc Proof:}}{\mbox{}\hfill 
\rule[-.2ex]{.25em}{1.8ex}}

\begin{document}

\title[Simple $C^{*}$-algebras with locally finite decomposition rank]{{\sc Simple $C^{*}$-algebras with locally finite decomposition rank}}

\author{Wilhelm Winter}
\address{Mathematisches Institut der Universit\"at M\"unster\\
Einsteinstr.\ 62\\ D-48149 M\"unster}

\email{wwinter@math.uni-muenster.de}

\date{\today}
\thanks{{\it Supported by:} EU-Network Quantum Spaces - Noncommutative 
Geometry (Contract No. \\
\indent HPRN-CT-2002-00280) and Deutsche 
Forschungsgemeinschaft (SFB 478)}

\setcounter{section}{-1}

\begin{abstract}
We introduce the notion of locally finite decomposition rank, a structural property shared by many stably finite nuclear $C^{*}$-algebras. The concept is particularly relevant for Elliott's program to classify nuclear $C^{*}$-algebras by $K$-theory data. We study some of its properties and show that a simple unital $C^{*}$-algebra, which has locally finite decomposition rank, real rank zero and which absorbs the Jiang--Su algebra $\Zh$ tensorially, has tracial rank zero in the sense of Lin. As a consequence, any such $C^{*}$-algebra, if it additionally satisfies the Universal Coefficients Theorem, is approximately homogeneous of topological dimension at most 3. Our result in particular confirms the Elliott conjecture for the class of simple unital $\Zh$-stable ASH algebras with real rank zero. Moreover, it implies that simple unital $\Zh$-stable AH algebras with real rank zero not only have slow dimension growth in the ASH sense, but even in the AH sense.   
\end{abstract}

\maketitle

\section{Introduction}

This note is concerned with the stably finite real rank zero case of Elliott's program to classify nuclear $C^{*}$-algebras by $K$-theory data; see \cite{R1} for an introduction to this subject. There is growing body of evidence that one can only expect $K$-theoretical classification results up to $\Zh$-stability, where $\Zh$ denotes the Jiang--Su algebra constructed in \cite{JS1} and a $C^{*}$-algebra $A$ is called $\Zh$-stable if it absorbs $\Zh$ tensorially. Salient results supporting this point of view can be found in \cite{GJS}, \cite{J}, \cite{PT}, \cite{R2}, \cite{T1}, \cite{T2}, \cite{TW1} and \cite{TW2}. \\
It is known that a $\Zh$-stable $C^{*}$-algebra $A$ behaves very well in many respects. In particular, it is either stably finite or purely infinite and, when exact, has nice comparison properties (cf.\ \cite{R2}). Moreover, $A$ has real rank zero if and only if the positive part of the $K_{0}$-group, $K_{0}(A)_{+}$, has dense image in the positive continuous affine functions on the tracial state space, $\Aff(T(A))_{+}$ (recall that $A$ has real rank zero if positive elements with finite spectrum are norm-dense in the set of all positive elements). \\ 
In \cite{W5} we confirmed the Elliott conjecture for the class of simple, separable, unital $C^{*}$-algebras which are $\Zh$-stable, have real rank zero and finite decomposition rank (to be explained below) and, additionally, satisfy the Universal Coefficients Theorem (UCT). In the present paper we generalize this result   to $C^{*}$-algebras which only have locally finite decomposition rank as opposed to finite decomposition rank. The difference might seem subtle at first glance, but we think that the generalization is substantial. The main point is that we use $\Zh$-stability  instead of  a condition like slow (or no) dimension growth  --  in our opinion this lends credibility to the point of view outlined above. 

Decomposition rank is a notion of covering dimension for nuclear $C^{*}$-algebras; it was introduced by E.\ Kirchberg and the author in \cite{KW}. Below we study a modified version of this concept: we say a $C^{*}$-algebra $A$ has locally finite decomposition rank if it can be exhausted by $C^{*}$-subalgebras each of which has finite decomposition rank. Note that we do not ask the decomposition ranks of the exhausting algebras to be globally bounded. Locally finite decomposition rank passes to quotients, inductive limits and to  hereditary subalgebras which are generated by projections; it implies nuclearity and quasidiagonality. Examples include all  separable approximately homogeneous (AH) $C^{*}$-algebras (in particular, all separable commutative $C^{*}$-algebras). In \cite{NgW}, P.\  W.\ Ng and the author have shown that separable approximately subhomogeneous (ASH) $C^{*}$-algebras also have locally finite decomposition rank. Clearly,  finite decomposition rank implies its local version.\\
We wish to emphazise that locally finite decomposition rank is a fairly mild condition on a stably finite nuclear $C^{*}$-algebra; it does not even exclude the known counterexamples to the Elliott conjecture in the stably finite case. In particular, it does not imply stable rank one, Blackadar's second fundamental comparability property or weak unperforation of the ordered $K_{0}$-group. These are all properties known to hold for nuclear stably finite $\Zh$-stable $C^{*}$-algebras by the results of \cite{R2}. In the classification results of \cite{El3}, \cite{EGL2} and \cite{W4}  (to mention but a few), they are guaranteed by conditions involving noncommutative  covering dimension, such as slow dimension growth or finite decomposition rank. In \cite{W5}, said properties were entailed by $\Zh$-stability, but, following the lines of \cite{W4}, we could as well have used our assumptions of finite decomposition rank and real rank zero  to obtain them -- a redundance which is removed in the present article. The main use of $\Zh$-stability in \cite{W5} was to get rid of a condition on the tracial state space still present in \cite{W4}. In fact, in the case of a unique tracial state the other hypotheses (real rank zero, $\Zh$-stability  and finite decomposition rank) can be considerably weakened as shown by N.\ Brown in \cite{B} and, more recently, by H.\ Lin in \cite{Li6}.  \\
Our main result generalizes Theorem 4.1 of \cite{W5}; it says that separable simple unital $\Zh$-stable $C^{*}$-algebras with locally finite decomposition rank and real rank zero have tracial rank zero. Using results of H.\ Lin, this confirms the Elliott conjecture for the class of such algebras which, additionally, satisfy the UCT. In particular, this applies to simple unital $\Zh$-stable ASH algebras with real rank zero. Thanks to earlier work of M.\ Dadarlat, G.\ Elliott, G.\ Gong and others, it then follows that such algebras are in fact AH of topological dimension at most 3 and that they have decomposition rank at most 2.

The paper is organized as follows: In Section 1 we introduce the concept of locally finite decomposition rank, study some of its properties and consider a number of examples. In Section 2 we state our main result and derive its corollaries. In the following section we outline our strategy for the proof of Theorem \ref{lfdrtr0} and describe the technical difficulties. Section 4 recalls some facts about order zero maps and $C^{*}$-algebras with real rank zero. Section 5 contains the key technical steps (Corollary \ref{excisible-appr} and Lemma \ref{excision}) for the proof of Theorem \ref{lfdrtr0}, which is completed in Section 6. \\

To try to prove a classification result using locally finite decomposition rank as opposed to finite decomposition rank was already suggested to me by Nate Brown several years ago, however, at that time I did not know how to use $\Zh$-stability to make such an attempt work. I would like to thank Nate as well as Ping Wong Ng, Mikael R{\o}rdam and Andrew Toms for many inspiring conversations on the classification program in general and and on $\Zh$-stability in particular.

\section{Locally finite decomposition rank}

Below we introduce the notion of locally finite decomposition rank, study some of its properties, compare it to the original decomposition rank and give a list of examples.

\bn
For convenience, we recall the following definition from \cite{KW}:

\label{d-dr} 
\begin{ndefn} (cf.\ \cite{KW}, Definitions 2.2 and 3.1) Let $A$ be a separable $C^*$-algebra.
\begin{itemize}
\item[(i)] A completely positive map $\varphi :  F \to A$ has   order zero, $\ord \varphi = 0$, if it preserves orthogonality, i.e., $\varphi(e) \varphi(f) = \varphi(f) \varphi(e) = 0$ for all $e,f \in F$ with $ef = fe = 0$.
\item[(ii)] A completely positive map $\varphi : F \to A$ ($F$ a finite-dimensional $C^{*}$-algebra) is $n$-decomposable, if there is a decomposition $F=F^{(0)} \oplus \ldots \oplus F^{(n)}$ such that the restriction of $\varphi$ to $F^{(i)}$ has   order zero for each $i \in \{0, \ldots, n\}$; we say $\varphi$ is $n$-decomposable with respect to $F=F^{(0)} \oplus \ldots \oplus F^{(n)}$.
\item[(iii)] $A$ has decomposition rank $n$, $\dr A = n$, if $n$ is the least integer such that the following holds: For any finite subset $\Gh  \subset A$ and $\varepsilon > 0$, there is a completely positive approximation $(F, \psi, \varphi)$ for $\Gh$ within $\varepsilon$ (i.e., $\psi:A \to F$ and $\varphi:F \to A$ are completely positive contractive and $\|\varphi \psi (b) - b\| < \varepsilon \; \forall \, b \in \Gh$) such that $\varphi$ is $n$-decomposable. If no such $n$ exists, we write $\dr A = \infty$.  
\end{itemize}
\end{ndefn}
\en

\bn
$C^{*}$-algebras with finite decomposition rank enjoy many nice properties (cf.\ \cite{KW}, \cite{W4}), but in some situations it would be desirable to have  a  condition which is fulfilled by a larger class of $C^{*}$-algebras, yet retains at least some of the nice structural properties implied by finite decomposition rank. There are several reasonable ways of weakening Definition \ref{d-dr}(iii). For example, one might ask the map $\varphi$  only to be completely positive contractive; this yields nothing but the completely positive approximation property, which is well-known to characterize nuclear $C^{*}$-algebras. An a priori less general version would be to ask the map $\varphi$ to have order zero on each of the summands of $F$; this definition does not rule out infinite $C^{*}$-algebras -- it might even be equivalent to the completely positive approximation property.\\
In these notes, we study a definition which does not entirely drop the decomposability condition of \ref{d-dr}(iii), but which also does not ask for a global bound on the decomposition constant:     

\begin{ndefn}
We say $A$ has locally finite decomposition rank, if, for any finite subset $\Gh \subset A$ and  $\varepsilon>0$, there is a $C^{*}$-subalgebra $B \subset A$ such that $\dr B$ is finite and  $\dist(b,B)<\varepsilon$ for all $b \in \Gh$.  
\end{ndefn}

Just like finite decomposition rank, this notion is a so-called \emph{local} property -- in fact, these two concepts may be thought of as local analogues of topologically finite-dimensional AH algebras and general AH algebras, respectively. We shall return to this point of view in \ref{lfdr-examples}.

\en

\bn
\label{lfdr-permanence}
\begin{nprop}
The property of having locally finite decomposition rank passes to inductive limits, quotients, tensor products and to  hereditary $C^{*}$-subalgebras generated by projections.
\end{nprop}

\begin{nproof}
The statements about limits, quotients and tensor products follows immediately from the respective statements for  decomposition rank, cf.\ \cite{W1}, Section 3, and \cite{KW}, 3.2. \\
Suppose $p$ is a projection in a $C^{*}$-algebra $A$ with locally finite decomposition rank. Let $\Gh \subset pAp$ be a  finite subset  and $\varepsilon>0$. We may assume that the elements of $\Gh$ are positive and normalized and that $p\in \Gh$. By assumption, for any $0<\delta<\varepsilon/3$ there is a $C^{*}$-subalgebra $B \subset A$ such that $\dr B < \infty$ and $\dist(b,B)<\delta \; \forall \, b \in \Gh$. But then it is straightforward to show that, if $\delta$ is chosen small enough, there is a partial isometry $s \in A$ such that $s^{*}s=p$, $q:=ss^{*} \in B$ and $\|s-p\|< \varepsilon/3$. Now $C:=s^{*}Bs$ is a $C^{*}$-subalgebra of $pAp$; since $s$ is a partial isometry, we have $C\cong qBq$. For any $b \in \Gh$, we have 
\[
\dist(b,C) = \dist(sbs^{*},qBq) \le \dist(b,qBq) + 2 \varepsilon/3 \le \dist(b,B) + 2 \varepsilon/3 \le \varepsilon \, ;
\]
by \cite{KW}, Proposition 3.8, $\dr C = \dr(qBq) \le \dr B < \infty$. We have thus shown that $pAp$ has locally finite decomposition rank.
\end{nproof}
\en

\bn
\label{nuclear-qd}
\begin{nprop}
A separable $C^{*}$-algebra $A$ with locally finite decomposition rank is nuclear and strongly quasidiagonal (i.e., every representation of $A$ is quasidiagonal); in particular, $A$ is stably finite.
\end{nprop}

\begin{nproof}
Since $A$ is exhausted by $C^{*}$-algebras with the completely positive approximation property, $A$ also has this property and hence is nuclear. \\
By \cite{BK2}, Corollary 5.7, a separable nuclear $C^{*}$-algebra is strongly quasidiagonal iff every quotient is strong NF in the sense of \cite{BK1}. By Proposition \ref{lfdr-permanence}, locally finite decomposition rank passes to quotients, so it will suffice to show that locally finite decomposition rank implies being strong NF. From \cite{KW}, Theorem 5.3, we already know that finite decomposition rank implies strong NF, and since being strong NF is a local property (see \cite{BK2}, Proposition 4.1 and the remark thereafter), the assertion follows.
\end{nproof}
\en

\bn
\label{lfdr-examples}
\begin{nexamples}
It is trivial that finite decomposition rank implies locally finite decomposition rank, so all the examples of \cite{KW}, Section 4,  of \cite{W3}, Section 1, and of \cite{TW2} have this property; this list includes the examples covered by virtually all known classification results for simple stably finite nuclear $C^{*}$-algebras. For example, all AF algebras, irrational rotation algebras and the Jiang--Su algebra $\Zh$ have (locally) finite decomposition rank.\\
There is a slight ambiguity in the literature about how to define approximately (sub-)homogeneous $C^{*}$-algebras (cf.\ \cite{Bl2}). We shall use the following set of definitions: A $C^{*}$-algebra $A$ is homogeneous, if all its irreducible representations have the same dimension. $A$ is approximately homogeneous (AH), if it is an inductive limit of direct sums of homogeneous $C^{*}$-algebras. $A$ is subhomogeneous, if the dimensions of its irreducible representations have some finite upper bound, and $A$ is approximately subhomogeneous (ASH), if it is an inductive limit of subhomogeneous $C^{*}$-algebras.\\
By \cite{Bl2}, Proposition 2.2, any separable AH algebra $A$ can be written as an inductive limit of direct sums of homogeneous algebras $A_{i}$ each of which has finite topological dimension (hence finite decomposition rank). Therefore, any AH algebra has locally finite decomposition rank, regardless of whether it has no, slow or fast dimension growth. In particular, this holds for Villadsen's examples and for Toms' counterexamples to the Elliott conjecture (cf.\ \cite{T1}). \\
In \cite{NgW}, Ping Wong Ng and the author showed the respective statements for ASH algebras, i.e., any separable ASH algebra is an inductive limit  $A= \lim_{\to} A_{i}$ of ASH algebras with finite topological dimension -- in particular, it has locally finite decomposition rank. Note that, again, we do not require the numbers $\dr A_{i}$ to have a common upper bound or the inductive limit decomposition to have slow dimension growth. 
\end{nexamples}
\en

\section{The main result and its consequences}

\bn
\label{lfdrtr0}
The concept of tracial rank zero was introduced by Lin (cf.\ \cite{Li2}, \cite{Li3}) as a somewhat more axiomatic approach to the stably finite real rank zero case of the Elliott program. We shall not need the original definition here (cf.\ \cite{Li0}, Definition 3.6.2), but we will give an alternative characterization in the next section, where we also outline the proof of the theorem below (the actual proof will have to wait until Section 6). Our main result states that many simple real rank zero $C^{*}$-algebras indeed have tracial rank zero:

\begin{ntheorem}
Let $A$ be a separable simple and unital $C^{*}$-algebra which is $\Zh$-stable and has real rank zero and locally finite decomposition rank. Then, $A$ has tracial rank zero.
\end{ntheorem}
\en

In \cite{Li2}, Lin confirmed the Elliott conjecture for the class of simple $C^{*}$-algebras with tracial rank zero which satisfy the UCT. We now explain how Lin's classification theorem for tracially AF algebras and results of Elliott (in the ASH case) and Dadarlat, Elliott and Gong (in the AH case) may be used to  derive a number of corollaries of Theorem \ref{lfdrtr0}; this is done in essentially the same way as in \cite{W5}. Moreover, we partially answer two questions of \cite{TW2}.\\

\bn
\begin{ncor}
Let $A$ be a separable simple unital $C^{*}$-algebra such that $A \otimes \Zh$ has locally finite decomposition rank. Then, the following are equivalent:
\begin{itemize}
\item[(i)] $A \otimes \Zh$ has tracial rank zero 
\item[(ii)] $A\otimes \Zh$ has real rank zero 
\item[(iii)] the canonical image of $K_{0}(A\otimes \Zh)$ in $\Aff(T(A\otimes \Zh))$ is dense
\item[(iv)] the canonical image of $K_{0}(A\otimes \Zh)_{+}$ in $\Aff(T(A\otimes \Zh))_{+}$ is dense.  
\end{itemize}
\end{ncor}

\begin{nproof}
(i) implies (ii) by \cite{Li0}, Theorem 3.6.11, the converse follows from Theorem \ref{lfdrtr0} above. (ii) and (iii) are equivalent by Proposition 7.1 of \cite{R2}. Since $A \otimes \Zh$ is nuclear and stably finite by Proposition \ref{nuclear-qd}, $A \otimes \Zh$ satisfies Blackadar's second fundamental comparability property by \cite{R2}, Corollary 4.10, whence (iii) and (iv) are equivalent.  
\end{nproof}
\en

\bn
\label{lfdr-classification}
In \cite{Li2}, Lin has confirmed the Elliott conjecture for the class of simple unital tracially AF algebras which satisfy the UCT. As a consequence we have the following

\begin{ncor}
Let $A$ and $B$ be separable simple unital $C^{*}$-algebras with real rank zero and locally finite decomposition rank; suppose $A$ and $B$ satisfy the UCT and are $\Zh$-stable. Then, $A$ and $B$ are isomorphic iff their Elliott invariants are.
\end{ncor}
\en

\bn
\label{AH-cor}
Thanks to the known results about the range of the Elliott invariant in the nuclear stably finite case, we can say more about the structure of algebras as in the preceding corollaries:

\begin{ncor}
Let $A$ be a separable simple unital $C^{*}$-algebra; suppose $A \otimes \Zh$ has real rank zero and locally finite decomposition rank and satisfies the UCT. Then:
\begin{itemize}
\item[(i)] $A \otimes \Zh$ is AH of topological dimension at most 3.
\item[(ii)] $A \otimes \Zh$ is ASH of topological dimension at most 2.
\item[(iii)] $\dr (A \otimes \Zh)$ is at most 2.
\item[(iv)] $A \otimes \Zh$ is approximately divisible. 
\item[(v)] $A$ is $\Zh$-stable iff $A$ is approximately divisible. 
\end{itemize}
\end{ncor}

\begin{nproof}
(i), (ii) and (iii) follow from results of Dadarlat, Elliott and Gong as in \cite{W4}, Corollary 6.4. By \cite{EGL}, an AH algebra of bounded topological dimension is approximately divisible. Conversely, an approximately divisible $C^{*}$-algebra is $\Zh$-stable by \cite{TW2}.
\end{nproof}
\en

\bn
We mention the following special case of Corollary \ref{lfdr-classification} explicitly:

\begin{ncor}
The class of separable simple unital $\Zh$-stable ASH $C^{*}$-algebras with real rank zero satisfies the Elliott conjecture. 
\end{ncor}

\begin{nproof}
ASH $C^{*}$-algebras clearly satisfy the UCT; they have locally finite decomposition rank by \cite{NgW}. The result follows from \ref{lfdrtr0} and \cite{Li2}. 
\end{nproof}
\en

\bn
\begin{nremarks}
(i) Note that \ref{AH-cor}(i) and (v) partially answer Questions 3.2 and 3.3 of \cite{TW2}. \\
(ii) In the preceding corollaries, note that the assumptions ``$A\otimes\Zh$ has real rank zero'' and ``$A\otimes \Zh$ has locally finite decomposition rank'' in particular hold if $A$ has real rank zero or locally finite decomposition rank, respectively (cf.\ Theorem 7.2 of \cite{R2}, Theorem 2.3 of \cite{W5} and Proposition \ref{lfdr-permanence} above).
\end{nremarks}
\en

\section{The proof of the main result: an outline}

Since we only have a rather complicated proof of Theorem \ref{lfdrtr0}, we outline our strategy below.

\bn
\label{wu-tr0}
First, we recall the definition of  simple tracial rank zero $C^{*}$-algebras in the presence of small projections and comparability. This characterization is an immediate consequence of \cite{Li0}, Definition 3.6.2 (cf.\ also \cite{Li3}, Corollary 6.15); it will be more useful for our purposes than the original definition.

\begin{nprop}
Let $A$ be a separable simple and unital $C^{*}$-algebra which satisfies Blackadar's second fundamental comparability property and every nonzero hereditary subalgebra of which contains a nonzero projection. Then, $A$ has tracial rank zero if and only if the following holds:\\
For any finite subset $\Fh \subset A$ and $\varepsilon>0$ there is a finite-dimensional $C^{*}$-subalgebra $D\subset A$ such that
\begin{itemize}
\item[(i)] $\|[\be_{D},b]\|<\varepsilon \; \forall \, b \in \Fh$
\item[(ii)] $\dist(\be_{D}b\be_{D},D)<\varepsilon \; \forall \, b \in \Fh$
\item[(iii)] $\tau(\be_{A} - \be_{D})< \varepsilon \; \forall \, \tau \in T(A)$.  
\end{itemize}
\end{nprop}
\en

\bn
\label{outline}
A $C^{*}$-algebra $A$ as in Theorem \ref{lfdrtr0} satisfies the hypotheses of the preceding proposition by results of R{\o}rdam (\cite{R2}). Therefore, given $\Fh \subset A$ and $\varepsilon>0$, we have to find a  finite-dimensional $C^{*}$-subalgebra $D \subset A$ satisfying (i), (ii) and (iii) above. \\
Since $A$ has locally finite decomposition rank, we may assume the elements of $\Fh$ to lie in some (unital) $C^{*}$-subalgebra $B$ of $A$ such that $\dr B=n$ for some $n \in \N$. Now suppose $B \stackrel{\psi}{\to} F \stackrel{\varphi}{\to} B$ is an $n$-decomposable  c.p.\ approximation of $\Fh$ within some $\alpha>0$. Since $A$ has real rank zero, we may replace $\varphi:F \to B$ by a so-called discretely $n$-decomposable map $\tilde{\varphi}:\tilde{F} \to A$ (cf.\  \ref{discrete-order-zero} and  \ref{rr0dr} below); the point is that $\tilde{\varphi} \circ \psi$ still is a good approximation for $\Fh$, while the image of $\tilde{\varphi}$ consists of a sum of $n+1$ (not necessarily pairwise orthogonal) finite-dimensional $C^{*}$-algebras $\tilde{F}^{(0)},\ldots,\tilde{F}^{(n)}$. Similar as in Section 4 of \cite{W5}   (using \ref{tracial-division} below), one can then use $\Zh$-stability of $A$ to find \emph{pairwise orthogonal} $C^{*}$-subalgebras $\bar{F}^{(0)}, \ldots, \bar{F}^{(n)}$ of $A$ such that $\bar{F}^{(i)} \cong \tilde{F}^{(i)}$ for all $i$ and such that $D_{1}:= \bar{F}^{(0)}\oplus \ldots \oplus \bar{F}^{(n)}$ satisfies (i) and (ii) above (with $D_{1}$ in place of $D$), if only $\alpha$ was chosen small enough. \\
This construction will not force $\bar{F}$ to quite satisfy (iii) -- the method of \cite{W5}, Section 4, will only yield $\tau(\be_{D_{1}}) > \frac{1}{2 (n+1)} =:\mu \; \forall \, \tau \in T(A)$. However, we may try to repeat the above process with $C_{1}:=(\be_{A} - \be_{D_{1}})A (\be_{A} - \be_{D_{1}})$ in place of $A$ and $\Fh_{1}:= \{(\be_{A} - \be_{D_{1}}) a (\be_{A} - \be_{D_{1}}) \, | \, x \in \Fh\}$ in place of $\Fh$ to obtain a finite-dimensional $D_{2} \supset D_{1}$ which not only satisfies (i) and (ii), but also $\tau(\be_{D_{2}}) > \mu (1-\mu) \; \forall \, \tau \in T(A)$. Induction will then yield an increasing sequence $D_{1} \subset D_{2} \subset \ldots \subset A$ such that $\tau(\be_{D_{k}}) > \mu \sum_{i =0}^{k} (1-\mu)^{i}$; by the formula for the geometric series we have $\mu \sum_{i =0}^{\infty} (1-\mu)^{i} = 1$, whence $\tau(\be_{D_{K}}) > 1- \varepsilon$ for some large enough $K$. \\
If $A$ itself has finite decomposition rank $n$, then all this works and, in fact, was carried out in \cite{W5}. But in the case where $A$ only has locally finite decomposition rank,  there is a major problem with the induction process: Although the algebras $C_{k}:=(\be_{A} - \be_{D_{k}})A (\be_{A} - \be_{D_{k}})$ again satisfy the same hypotheses as $A$, we can only be sure to be able to approximate the elements of $\Fh_{k}$ by $m$-decomposable c.p.\ approximations for some $m \in \N$, but it may well happen that $m$ is much larger than $n$ -- and this would destroy the final geometric series argument.  The difficulty could be circumvented if the compression with  $(\be_{A} - \be_{D_{k}})$ was multiplicative on $B$, for then the image of $B$ in $C_{k}$ again had decomposition rank $n$ and we could proceed as before by approximating the elements of $\Fh_{k}$ with $n$-decomposable c.p.\ approximations. Of course, in general compression with $(\be_{A} - \be_{D_{k}})$ will not be multiplicative -- but with the help of (i) and (ii) above  (with improved approximation constants) we can assume it to be \emph{almost} multiplicative with respect to some tolerance and some finite subset (which includes $\Fh_{k}$ and $\varphi(F)$). This will still be enough to obtain an $n$- (as opposed to $m$-) decomposable c.p.\ approximation of $\Fh_{k}$, and it will allow our induction process to work. The latter assertion is (roughly speaking) the content of our technical key results, \ref{excisible-appr} and \ref{excision}, the proof of which is the objective of Section 5. \\
What makes this procedure so complicated is the necessity to carefully keep track of the approximation constants chosen along the way. In fact, given $\Fh$ and $\varepsilon$, we first chose $B$ and, at the same time, obtain $n$. This $n$ determines how many induction steps will be needed ($\mu \sum_{i =0}^K (1-\mu)^{i} $ has to be larger than $1- \varepsilon$, and $\mu$ depends on $n$). Next we choose $\alpha$ and the c.p.\ approximation $(F,\psi,\varphi)$. The number $\alpha$ has to be so small that, even after $K$ induction steps,  the algebra $D_{K}$ still satisfies (i) and (ii) of Proposition \ref{wu-tr0} (this is where Lemma \ref{excision} enters).  Only now we can let the induction process start, i.e.,  carry out the actual construction of the $D_{k}$ for $k=1, \ldots,K$. These last steps will complete the proof of Theorem \ref{lfdrtr0} and are the content of Section 6.          
\en

\bn
One might ask whether some of the technicalities outlined above could be avoided by using ultraproduct techniques. Such an approach could in fact help to replace the above mentioned compression with $(\be_{A} - \be_{D_{k}})$ by an honestly multiplicative map into the ultraproduct $A_{\omega}$ ($\omega$ being some free ultrafilter  on $\N$). However, it does not seem to be possible to carry out the whole induction procedure of \ref{outline} just in the ultraproduct -- one would rather have to lift the multiplicative map into $A_{\omega}$ to a sequence of almost multiplicative maps into $A$, and this would leave us in essentially the same situation as before, so the technical advantages of employing ultraproducts seem to be rather moderate. Nontheless, such an approach will be used in \cite{NgW2} to prove a result related to our Theorem \ref{lfdrtr0}.  
\en

\section{Order zero maps}

In this section we recall some facts about $n$-decomposable maps into $C^{*}$-algebras of real rank zero.

\bn
\label{discrete-order-zero}
Recall from \cite{W4}, Definition 2.2(i), that a completely positive map 
\[
\varphi: F=M_{r_{1}} \oplus \ldots \oplus M_{r_{s}} \to A 
\] 
is a discrete order zero map, if $\ord \varphi = 0$ and each $\varphi(\be_{M_{r_{i}}})$, $i=1, \ldots, s$, is a multiple of a projection. \\
Let  $\tilde{F}$ be another finite-dimensional $C^{*}$-algebra. We say an embedding $\iota:F \to \tilde{F}$ is centered, if there are $m_{1}, \ldots,m_{s} \in \N$ such that $\tilde{F}\cong \bigoplus_{i=1}^{s} \C^{m_{i}} \otimes M_{r_{i}}$ and, under this identification, 
\[
\iota = \bigoplus_{i=1}^{s} \be_{\C^{m_{i}}} \otimes \id_{M_{r_{i}}} \, .
\]
This is equivalent to saying that the commutant of $\iota(F)$ within $\tilde{F}$ coincides with the center of $\tilde{F}$.
\en

\bn
By \cite{W4}, Lemma 2.4 (and its proof), any order zero map into a real rank zero $C^{*}$-algebra $A$ can be approximated by a composition of a centered embedding with a discrete order zero map:

\begin{nlemma}
Let $A$ and $F$ be $C^{*}$-algebras, $A$ with real rank zero and $F$ finite-dimensional. Suppose $\varphi:F \to A$ is completely positive contractive with   order zero and let $\delta>0$ be given. Then there are a centered unital embedding $\iota:F \to \tilde{F}$ of $F$ into some finite-dimensional $C^{*}$-algebra $\tilde{F}$ and a discrete order zero map $\tilde{\varphi}:\tilde{F} \to A$ such that $\tilde{\varphi}(\be_{\tilde{F}})\le \varphi(\be_{F})$ and $\|\varphi(x) - \tilde{\varphi} \circ \iota(x)\|< \delta \cdot \|x\|$ for all $0 \neq x \in F$. 
\end{nlemma}
\en

\bn
\label{rr0dr}
The preceding Lemma carries over to $n$-decomposable maps, as the next proposition shows. First, we need some notation: Let $A$ and $F$ be $C^{*}$-algebras with $F$ finite-dimensional, and let $\varphi:F \to A$ be a c.p.\ map. Following \cite{W4}, Definition 2.2(ii), we say $\varphi$ is discretely $n$-decomposable, if $F$ can be written as $F = F^{(0)} \oplus \ldots \oplus F^{(n)}$ with $\varphi|_{F^{(j)}}$ being a discrete order zero map for $j=0, \ldots, n$.  

\begin{nprop}  Let $A$ and $F$ be  $C^{*}$-algebras, $F=M_{r_{1}} \oplus \ldots \oplus M_{r_{s}}$ finite-dimensional and $A$ with real rank zero. Let $\varphi:F \to A$ be an $n$-decomposable c.p.c.\ map.\\
Then, for any $\beta>0$, there are a centered unital embedding $\iota:F \to \tilde{F}$ into some finite-dimensional $C^{*}$-algebra $\tilde{F}$ and a discretely $n$-decomposable c.p.c.\ map $\tilde{\varphi}:\tilde{F} \to A$ such that $\tilde{\varphi}\circ \iota(\be_{M_{r_{i}}})\le \varphi(\be_{M_{r_{i}}})$ and $\|\varphi(x) - \tilde{\varphi} \circ \iota(x)\|< \beta \cdot \|x\|$ for all $i =1, \ldots, s$ and $0 \neq x \in F$.\\
If $\varphi$ is $n$-decomposable with respect to the decomposition $F=F^{(0)} \oplus \ldots \oplus F^{(n)}$, then $\tilde{\varphi}$ may be chosen to be $n$-decomposable with respect to the decomposition $\tilde{F}=\tilde{F}^{(0)} \oplus \ldots \oplus \tilde{F}^{(n)}$, where $\tilde{F}^{(j)}={\iota}(F^{(j)})$, $j=0, \ldots,n$. 
\end{nprop}

\begin{nproof}
Apply Lemma \ref{discrete-order-zero} with $\delta:=\frac{\beta}{n+1}$ to each of the maps $\varphi|_{M_{r_{i}}}$ to obtain discrete order zero maps $\tilde{\varphi}_{i}$, $i=1, \ldots, s$. The $\tilde{\varphi}_{i}$ will add up to a discretely $n$-decomposable map $\tilde{\varphi}$ with the desired properties; cf.\ also the proof of \cite{W4}, Proposition 2.5. 
\end{nproof}
\en

\bn
\label{multiplicative-domain}
We shall have use for the following consequence of Stinespring's theorem, which  is a standard tool to analyze completely positive approximations of nuclear $C^{*}$-algebras. See \cite{KW}, Lemma 3.5, for a proof.
 
\begin{nlemma}
Let $A$ and $F$ be $C^{*}$-algebras, $b \in A$ a normalized positive element and $\eta>0$. If $A \stackrel{\psi}{\longrightarrow} F \stackrel{\varphi}{\longrightarrow} A$ are completely positive contractive maps satisfying 
\[
\|\varphi \psi(b) - b\|, \, \|\varphi \psi(b^{2}) - b^{2}\| < \eta \, ,
\]
then, for any $0 \neq x \in F_{+}$, 
\[
\|\varphi(\psi(b)x)- \varphi \psi(b) \varphi(x)\|< 2 \eta^{\halb} \|x\| \, .
\]
\end{nlemma}
\en

\bn
The proof of Theorem \ref{lfdrtr0} becomes considerably easier in the case of finitely (or countably) many tracial states. The following lemma  (2.4 from \cite{W5}) will be used to avoid this assumption:

\label{tracial-division}
\begin{nlemma}
For any $n \in \N$ and $0< \mu <  1/2(n+1)$ there is a completely positive contractive   order zero map $\varrho: \C^{n+1} \to \Zh$ such that $\bar{\tau}(\varrho(e_{i})) > \mu$ for $i=1, \ldots, n+1$, where  the $e_{i}$ denote the canonical generators of $\C^{n+1}$ and $\bar{\tau}$ is the unique tracial state on $\Zh$. 
\end{nlemma}
\en

\section{Excising almost central subalgebras}

This section contains the technical key steps for the proof of Theorem \ref{lfdrtr0}, namely Corollary \ref{excisible-appr} and Lemma \ref{excision}. First, we need some preparation.

\bn
\label{polynomial-appr}
\begin{nprop}
For any $\delta>0$ and $f ,g\in \Ch_{0}((0,1])$ there is $0 < \beta < \delta$ such that the following holds: If $0\le a,b \le \be$ are elements in some $C^{*}$-algebra which satisfy $\|a-b\|<\beta$ (or $\|[a,b]\|< \beta$, respectively),  then $\|f(a)-f(b)\|< \delta$ (or $\|[g(a),f(b)]\|< \delta$, respectively). 
\end{nprop}

\begin{nproof}
The assertions are obvious if $f$ and $g$  are polynomials. By the Stone--Weierstrass Theorem any function in $\Ch_{0}((0,1])$ is a uniform limit of polynomials, from which the statements follow immediately. 
\end{nproof}
\en

\bn
\label{almost-hereditary}
\begin{nprop}
Let $0 \le a, b \le \be_{A}$ be positive elements of a unital $C^{*}$-algebra $A$ and let $\varepsilon>0$ be given.  If $a \le b + \varepsilon \cdot \be_{A}$, then $\dist (a, \overline{bAb}) \le 3 \cdot \varepsilon^{\halb}$.
\end{nprop}

\begin{nproof}
Let $(u_{n})_{n \in \N} \subset \overline{bAb}$ be an approximate unit of $\overline{bAb}$;  assume that $0 \le u_{n} \le \be_{A}$. We then have
\begin{eqnarray*}
\|a - u_{n} a u_{n}\| & \le & \|(\be_{A}-u_{n}) a u_{n}\| + \|u_{n} a (\be_{A} - u_{n})\| + \|(\be_{A}-u_{n}) a (\be_{A}-u_{n})\| \\
& \le & 2 \|(\be_{A}-u_{n}) a u_{n}^{2}a (\be_{A}- u_{n})\|^{\halb} + \|(\be_{A} - u_{n}) a (\be_{A} - u_{n})\| \\
& \le & 3 \|(\be_{A} - u_{n}) a (\be_{A} - u_{n})\|^{\halb} \\
& \le & 3 (\|(\be_{A} - u_{n}) b (\be_{A} - u_{n})\| + \varepsilon)^{\halb} \, ,
\end{eqnarray*}
from which follows that, for any $\delta>0$, there is $n \in \N$ such that 
\[
\|a - u_{n} a u_{n}\|< 3 (\delta + \varepsilon)^{\halb} \, .
\]
Since $u_{n}a u_{n} \in \overline{bAb}$ and $\delta$ is arbitrary, the assertion follows.
\end{nproof}
\en

\bn
\label{functions}
\begin{nnotation}
For $0 < \alpha <\beta < 1$ we define continuous functions 
\[
g_{\alpha,\beta},h_{\alpha,\beta}:[0,1] \to \R
\]
by
\[
g_{\alpha,\beta}(t) := \left\{
\begin{array}{ll}
0, & 0\le t \le {\alpha} \\
1, & \beta \le t \le 1 \\
\mbox{linear,} & \mbox{else}
\end{array}
\right.
\]
and
\[
h_{\alpha,\beta}(t) := \left\{
\begin{array}{ll}
0, & 0\le t \le {\alpha} \\
t^{-1}, & {\beta} \le t \le 1 \\
\mbox{linear,} & \mbox{else} \, .
\end{array}
\right.
\]
The subset of positive elements of norm at most one in a $C^{*}$-algebra $B$ will be denoted by $\Bh_{1}(B_{+})$.
\end{nnotation}
\en

\bn
\label{p-kappa}
\begin{nlemma}
Let $A$ be a unital $C^{*}$-algebra and $B \subset A$ a unital $C^{*}$-subalgebra. Furthermore, let $\Gh \subset \Bh_{1}(B_{+})$ be a compact subset containing $\be_{A}$ and let $n \in \N$ and $0 < \zeta < 1/19$ be given. Then, there is $\zeta'>0$ such that the following holds:\\
If $(F, \psi, \varphi)$ is an $n$-decomposable c.p.\ approximation of $B$ such that
\[
\|\varphi \psi(b) - b \| < \frac{\zeta^{6}}{(n+1)^{2}}\; \forall \, b \in \bar{\Gh}:=\Gh \cup \{a^{2} \, | \, a \in \Gh\}
\]
and if $F_{1}, \ldots ,F_{s}$ are the matrix blocks of $F$ and $p_{1}, \ldots, p_{s} \in A$ are pairwise orthogonal projections satisfying 
\begin{equation}
\label{34}
\|[p_{i},\varphi( \be_{F_{i}}x)]\| < \zeta'  \|x\| \; \forall \, 0 \neq x \in F_{+} \, ,
\end{equation}
\begin{equation}
\label{26}
\|p_{i} g_{\frac{\zeta}{2},\zeta}(\varphi(\be_{F_{i}})) - p_{i}\| < \zeta'
\end{equation}
and 
\[
\dist(p_{i} , \overline{\varphi(\be_{F_{i}}) A \varphi(\be_{F_{i}})}) < \frac{\zeta'}{s}
\]
for $i=1, \ldots , s$, then $p:= \sum_{i=1}^{s} p_{i}$ satisfies
\[
\|[p,b]\| <  \zeta 
\]
for all $b \in \Gh$. 
\end{nlemma}

\begin{nproof}
Consider $h_{\frac{\zeta}{4},\frac{\zeta}{2}} \in \Ch_{0}((0,1])$ and note that 
\begin{equation}
\label{27}
\id_{[0,1]} \cdot g_{\frac{\zeta}{2},\zeta} \cdot h_{\frac{\zeta}{4},\frac{\zeta}{2}} = g_{\frac{\zeta}{2},\zeta} 
\end{equation} 
and that
\begin{equation}
\label{31}
\|h_{\frac{\zeta}{4},\frac{\zeta}{2}}\| = \frac{2}{\zeta} \, .
\end{equation}

By Proposition \ref{polynomial-appr}, there is $\zeta'$ such that the following holds: If $0 \le a,b \le \be_{A}$ are elements of $A$ with $\|[a,b]\|< \zeta'$, then 
\begin{equation}
\label{35} \|[a, (g_{\frac{\zeta}{2},\zeta} \cdot h_{\frac{\zeta}{4},\frac{\zeta}{2}})(b)]\| < \frac{1}{19(n+1)} \zeta \, .
\end{equation}
We may assume that 
\begin{equation}
\label{36}
\zeta' < \frac{1}{38(n+1)} \zeta^{2} \, .
\end{equation} 

Now suppose that  $(F=F_{1}\oplus \ldots \oplus F_{s},\psi,\varphi)$ is a c.p.\ approximation and $p_{1}, \ldots, p_{s}\in A$ are projections as in the statement of the proposition.  Let $\varphi$ be $n$-decomposable with respect to the decomposition $F=(\bigoplus_{i \in I_{0}} F_{i}) \oplus \ldots \oplus (\bigoplus_{i \in I_{n}} F_{i})$, where $\{1, \ldots,s\} = \coprod_{j=0}^{n} I_{j}$; in particular, this means that, for all $j \in \{0, \ldots,n\}$, 
\begin{equation}
\label{29}
\varphi(\be_{F_{i}}) \perp \varphi(\be_{F_{i'}}) \mbox{ if } i \neq i' \in I_{j}
\end{equation}
and
\begin{equation}
\label{32}
[\varphi(\be_{F_{i}}), \varphi(\be_{F_{i}} x)] = 0 \; \forall \, i \in 1, \ldots, s, \, x \in F \, .
\end{equation}
By Lemma \ref{multiplicative-domain} and our assumption on $(F,\psi,\varphi)$ we have 
\begin{equation}
\label{30}
\| (\sum_{i \in I_{j}} \varphi(\be_{F_{i}}) ) \varphi\psi(b) - \varphi( \sum_{i \in I_{j}} \be_{F_{i}} \psi(b)) \| < 2 \cdot \frac{\zeta^{3}}{n+1} 
\end{equation}
for each $j \in \{0, \ldots,n\}$ and $b \in \Gh$. Since $\dist(p_{i},\overline{\varphi(\be_{F_{i}})A\varphi(\be_{F_{i}})}) < \frac{\zeta'}{s}$ for each $i$, there are positive normalized elements 
\begin{equation}
\label{33}
d_{i} \in C^{*}(\varphi(\be_{F_{i}}))
\end{equation}
such that 
\begin{equation}
\label{25}
\|p_{i}-d_{i}p_{i}d_{i}\|<  \frac{\zeta'}{s} \; \forall \, i \, . 
\end{equation}
For each $j \in \{0, \ldots , n\}$ we obtain
\begin{eqnarray}
\label{28}
\lefteqn{
\|\sum_{i \in I_{j}} p_{i} - \sum_{i \in I_{j}} d_{i} p_{i} d_{i} (g_{\frac{\zeta}{2},\zeta} \cdot h_{\frac{\zeta}{4},\frac{\zeta}{2}})(\varphi(\be_{F_{i}}))\varphi(\be_{F_{i}})\|
} \nonumber \\
& \stackrel{(\ref{25})}{\le} & \|\sum_{i \in I_{j}} d_{i}p_{i}d_{i} - \sum_{i \in I_{j}} d_{i} p_{i} d_{i} (g_{\frac{\zeta}{2},\zeta} \cdot h_{\frac{\zeta}{4},\frac{\zeta}{2}})(\varphi(\be_{F_{i}}))\varphi(\be_{F_{i}})\| + s \cdot \frac{\zeta'}{s} \nonumber \\
& \stackrel{(\ref{29},\ref{32},\ref{33})}{\le} & \max_{i \in I_{j}} \| d_{i} ( p_{i} -   p_{i}  (g_{\frac{\zeta}{2},\zeta} \cdot h_{\frac{\zeta}{4},\frac{\zeta}{2}})(\varphi(\be_{F_{i}}))\varphi(\be_{F_{i}}) d_{i} \| + \zeta' \nonumber \\
& \stackrel{(\ref{26},\ref{27})}{\le} & 2 \cdot \zeta' \, .
\end{eqnarray}
We now compute for any $b \in \Gh$
\begin{eqnarray*}
\lefteqn{
\|[(\sum_{i \in I_{j}} p_{i}), \varphi \psi(b)]\|
}\\
& \stackrel{(\ref{28})}{\le} & \|\sum_{i \in I_{j}} d_{i}p_{i}d_{i} (g_{\frac{\zeta}{2},\zeta} \cdot h_{\frac{\zeta}{4},\frac{\zeta}{2}})(\varphi(\be_{F_{i}})) \varphi(\be_{F_{i}}) \varphi \psi(b) \\
& & -  \varphi \psi(b)  \sum_{i \in I_{j}} \varphi(\be_{F_{i}})   (g_{\frac{\zeta}{2},\zeta} \cdot h_{\frac{\zeta}{4},\frac{\zeta}{2}})(\varphi(\be_{F_{i}})) d_{i}p_{i}d_{i}\| \\
& & + 4 \zeta' \\
& \stackrel{(\ref{29})}{=} & \|\sum_{i \in I_{j}} d_{i}p_{i}d_{i} (g_{\frac{\zeta}{2},\zeta} \cdot h_{\frac{\zeta}{4},\frac{\zeta}{2}})(\varphi(\be_{F_{i}})) \sum_{i \in I_{j}} \varphi(\be_{F_{i}}) \varphi \psi(b) \\
& & -  \varphi \psi(b)  \sum_{i \in I_{j}} \varphi(\be_{F_{i}})   \sum_{i \in I_{j}}(g_{\frac{\zeta}{2},\zeta} \cdot h_{\frac{\zeta}{4},\frac{\zeta}{2}})(\varphi(\be_{F_{i}})) d_{i}p_{i}d_{i}\| \\
& & + 4 \zeta' \\
& \stackrel{(\ref{30},\ref{31})}{\le} & \|(\sum_{i \in I_{j}} d_{i}p_{i}d_{i} (g_{\frac{\zeta}{2},\zeta} \cdot h_{\frac{\zeta}{4},\frac{\zeta}{2}})(\varphi(\be_{F_{i}}))) \varphi(\sum_{i \in I_{j}} \be_{F_{i}} \psi(b)) \\
& & -  \varphi (\sum_{i \in I_{j}} \psi(b) \be_{F_{i}})   \sum_{i \in I_{j}}(g_{\frac{\zeta}{2},\zeta} \cdot h_{\frac{\zeta}{4},\frac{\zeta}{2}})(\varphi(\be_{F_{i}})) d_{i}p_{i}d_{i}\| \\
& & + 4 \zeta' + 2 \cdot \frac{2}{\zeta} \cdot 2 \cdot \frac{\zeta^{3}}{n+1}\\ 
& \stackrel{(\ref{29})}{=} & \|\sum_{i \in I_{j}}[d_{i}p_{i}d_{i}, (g_{\frac{\zeta}{2},\zeta} \cdot h_{\frac{\zeta}{4},\frac{\zeta}{2}})(\varphi(\be_{F_{i}}))\varphi(\be_{F_{i}}\psi(b))]\| \\
&& + 4 \zeta' + 8 \cdot \frac{\zeta^{2}}{n+1}  \\ 
& \stackrel{(\ref{32},\ref{33})}{=} & \|\sum_{i\in I_{j}} d_{i} [p_{i}, (g_{\frac{\zeta}{2},\zeta} \cdot h_{\frac{\zeta}{4},\frac{\zeta}{2}})(\varphi(\be_{F_{i}}))\varphi(\be_{F_{i}}\psi(b))] d_{i}\| \\
&& + 4 \zeta' + 8 \cdot \frac{\zeta^{2}}{n+1} \\
& \stackrel{(\ref{29},\ref{33})}{=} & \max_{i \in I_{j}} \| d_{i} [p_{i}, (g_{\frac{\zeta}{2},\zeta} \cdot h_{\frac{\zeta}{4},\frac{\zeta}{2}})(\varphi(\be_{F_{i}}))\varphi(\be_{F_{i}}\psi(b))] d_{i}\| \\
&& + 4 \zeta' + 8 \cdot \frac{\zeta^{2}}{n+1} \\
& \le & \max_{i \in I_{j}} \|  [p_{i}, (g_{\frac{\zeta}{2},\zeta} \cdot h_{\frac{\zeta}{4},\frac{\zeta}{2}})(\varphi(\be_{F_{i}}))\varphi(\be_{F_{i}}\psi(b))] \| \\
&& + 4 \zeta' + 8 \cdot \frac{\zeta^{2}}{n+1} \\
& \stackrel{(\ref{34},\ref{35},\ref{31})}{\le} & \frac{\zeta}{19(n+1)} + \zeta' \cdot \frac{2}{\zeta} + 4 \zeta' + 8 \cdot \frac{\zeta^{2}}{n+1}  \, .
\end{eqnarray*}
As a consequence, we obtain 
\begin{eqnarray*}
\|[p,b]\| &\le &  \|[p, \varphi \psi(b)]\| + 2 \cdot \frac{\zeta^{6}}{(n+1)^{2}}\\
& \le & (n+1) \left(\frac{\zeta}{19(n+1)} + \zeta' \cdot \frac{2}{\zeta} + 4 \zeta' + 8 \cdot \frac{\zeta^{2}}{n+1}  \right) + 2 \cdot \frac{\zeta^{6}}{(n+1)^{2}}\\
& < & (n+1) \left(\frac{\zeta}{19(n+1)} + \zeta' \cdot \frac{2}{\zeta} + 4 \zeta' + 8 \cdot \frac{\zeta^{2}}{n+1}  + 2 \cdot \frac{\zeta^{6}}{n+1} \right)\\
& \stackrel{(\ref{36})}{<} & \zeta  
\end{eqnarray*}
for all $b \in \Gh$.
\end{nproof}
\en

\bn
\label{excisible-appr}
For convenience, we note the following corollary explicitly:

\begin{ncor}
Let $A$ be a unital $C^{*}$-algebra and $B \subset A$ a unital $C^{*}$-subalgebra with $\dr B = n < \infty$. \\
For any compact subset $\Gh \subset \Bh_{1}(B_{+})$, and $0 < \eta <1/19$ there are an $n$-decomposable c.p.\ approximation $(F, \psi, \varphi)$ of $B$ and $\delta>0$ such that the following hold:
\begin{itemize}
\item[a)] $\|\varphi \psi(b) - b\| < \frac{\eta^{6}}{(n+1)^{2}} \; \forall \, b \in \bar{\Gh}:= \Gh \cup \{a^{2} \, | \, a \in \Gh\}$
\item[b)] If $F_{1},  \ldots, F_{s}$ are the matrix blocks of $F$ and $p_{1}, \ldots,p_{s} \in A$ are pairwise orthogonal projections satisfying 
\[
\|[p_{i}, \varphi(\be_{F_{i}} x)]\| < \delta \|x\| \; \forall \, 0 \neq x \in F_{+}  \, ,
\]
\[
\|p_{i} g_{\frac{\eta}{2},\eta}(\varphi(\be_{F_{i}})) - p_{i}\| < \delta
\]
and
\[
\dist(p_{i} , \overline{\varphi(\be_{F_{i}})A \varphi(\be_{F_{i}})} ) < \frac{\delta}{s}
\]
for $i=1, \ldots, s$, then $p:= \sum_{i=1}^{s}p_{i}$ satisfies 
\[
\|[p,b]\| < \eta \; \forall \, b \in \Gh \, .
\]
\end{itemize}
\end{ncor}
\en

\bn
\label{excision}
\begin{nlemma}
Let $A$ be a separable simple and unital $\Zh$-stable $C^{*}$-algebra with real rank zero and let $B \subset A$ be a unital $C^{*}$-subalgebra. Let $(F,\psi,\varphi)$ be a c.p.\ approximation of $B$ and suppose $\varphi$ is $n$-decomposable for some $n\in \N$. Let $\mu$ and $\eta$ be positive numbers such that 
\begin{equation}
\label{62}
0<\mu<\frac{1}{2(n+1)}, \, \eta<\frac{1}{48} \mbox{ and }  \eta < \frac{1}{10} \left(\frac{1}{2(n+1)}-\mu \right) \, .
\end{equation}
Furthermore, let $\Gh \subset \Bh_{1}(B_{+})$ be a compact subset containing $\be_{A}$ and satisfying
\begin{equation}
\label{51}
\|\varphi \psi(b)-b\|<\frac{\eta^{6}}{(n+1)^{2}} \; \forall \,  b \in \bar{\Gh}:= \Gh \cup \{a^{2} \, | \, a \in \Gh\} \, .
\end{equation}
Then, for any $0< \delta < \halb$ there is $\gamma>0$ such that the following holds:\\
If there is a projection $q \in A$ such that 
\[
\|[q,\varphi(x)]\| < \gamma \|x\| \; \forall \, 0 \neq x \in F_{+} \, ,
\]
then there is a finite-dimensional $C^{*}$-subalgebra $C \subset (\be_{A}-q)A(\be_{A}-q) \subset A$ such that
\begin{itemize}
\item[(i)] $\dist(\be_{C} b \be_{C},C) < \eta \; \forall \, b \in \Gh$
\item[(ii)] $\tau(\be_{C}) \ge \mu \cdot \tau(\be_{A}-q) \; \forall \, \tau \in T(A)$
\item[(iii)]  if $F_{1}, \ldots, F_{s}$ are the matrix blocks of $F$, then $\be_{C}$ can be written as a sum of $s$ pairwise orthogonal projections $p_{1}, \ldots,p_{s} \in C$ satisfying
\[
\|[p_{i},\varphi(\be_{F_{i}}x)]\| < \delta \|x\| \; \forall \, 0 \neq x \in F_{+} \, ,  
\]
\[
\|p_{i} g_{\frac{\eta}{2},{\eta}}(\varphi(\be_{F_{i}})) - p_i\|<\delta
\]
and
\[
\dist (p_{i},  \overline{\varphi(\be_{F_{i}})A\varphi(\be_{F_{i}})} ) < \frac{\delta}{s}
\]
for $i=1, \ldots,s$.
\end{itemize}
\end{nlemma}

\begin{nproof}
For convenience, we define $g,h\in \Ch_{0}((0,1])$ by 
\begin{equation}
\label{44}
g:=g_{\frac{\eta}{2},\eta} \mbox{ and } h:= h_{\frac{\eta}{4},\frac{\eta}{2}} \cdot g_{\frac{\eta}{2},\eta} \, .
\end{equation}
We clearly have 
\begin{equation} 
\label{41}
t \cdot h(t) = g(t) \; \forall \, t \in [0,1] \mbox{ and } \|h\|=\frac{1}{\eta} \, .
\end{equation}

Given  $\delta $, use Proposition \ref{polynomial-appr} to choose $\beta>0$ such that if $a$ and $b$ are positive elements in some $C^{*}$-algebra which have norm at most one and satisfy $\|a-b\|<\beta$, then  
\begin{equation}
\label{38}
\|g(a) - g(b)\| < \frac{\delta^{2}}{12} \, .
\end{equation}
We may also assume that 
\begin{equation}
\label{48}
\beta< \frac{\delta^{2}\eta^{6}}{12s^{2}(n+1)^{2}} \, .
\end{equation}
By \cite{KW}, Remark 2.4 and Proposition 2.5, the relations defining $n$-decomposability are weakly stable in the sense of \cite{Lo}, Definition 4.1.1. This implies that there is $\gamma>0$ such that the following holds:\\
If there is a projection $q \in A$ such that 
\begin{equation}
\label{37}
\|[q,\varphi(x)]\| < \gamma \|x\| \; \forall \, 0 \neq x \in F_{+} \, ,
\end{equation}
then there are $n$-decomposable c.p.\ maps 
\begin{equation}
\label{61}
\varphi' : F \to (\be_{A} - q) A (\be_{A} - q)
\end{equation}
and 
\begin{equation}
\label{55}
\varphi^{\times}: F \to qAq
\end{equation}
such that
\begin{equation}
\label{13}
\|\varphi'(x) - (\be_{A} - q)\varphi(x) (\be_{A} - q) \| < \beta \|x\| \; \forall \, 0 \neq x \in F_{+} \, ,
\end{equation} 
\[
\|\varphi^{\times}(x) - q \varphi(x) q\| < \beta \|x\| \; \forall \, 0 \neq x \in F_{+} 
\]
and
\begin{equation}
\label{54}
\|\varphi'(x) +\varphi^{\times}(x) - \varphi(x)\|<\beta \|x\| \; \forall \, 0 \neq x \in F_{+} \, .
\end{equation}
(In other words, the c.p.\ maps $q\varphi(\, . \,)q$ and $(\be_{A}-q)\varphi(\, . \,)(\be_{A}-q)$ are `almost' $n$-decomposable whatever particular $q \in A$ we choose, if only (\ref{37}) is satisfied.) By making $\gamma$ smaller, if necessary, and using Proposition \ref{polynomial-appr} and the fact that 
\[
g(\varphi'(\be_{F_{i}})) \le g(\varphi'(\be_{F_{i}})) +g(\varphi^{\times}(\be_{F_{i}})) = g(\varphi'(\be_{F_{i}}) + \varphi^{\times}(\be_{F_{i}})) \, ,
\]
we may even assume that 
\begin{equation}
\label{71}
g(\varphi'(\be_{F_{i}})) \le g(\varphi(\be_{F_{i}})) + \frac{\delta^{2}}{12}
\end{equation}
for all $i \in \{1, \ldots, s\}$. We may further assume that 
\begin{equation}
\label{50}
\gamma < \frac{\delta^{2} \eta^{2}}{16 s^{2} (n+1)^{2}} \, .
\end{equation}

So, let $q \in A$ as above be given and suppose we have chosen $\varphi'$ and $\varphi^{\times}$. From Proposition \ref{rr0dr} and the choice of $\beta$ we obtain a finite-dimensional $C^{*}$-algebra $\bar{F}$ with a unital  centered embedding $\bar{\iota}: F \to \bar{F}$ and a discretely $n$-decomposable c.p.c.\ map 
\[
\varphi'': \bar{F} \to (\be_{A} - q) A (\be_{A} - q)
\]
such that
\begin{equation}
\label{4}
\|\varphi'' \circ \bar{\iota}(x) - \varphi'(x) \| < \beta \|x\| \; \forall \, 0 \neq x \in F_{+} 
\end{equation}
and
\begin{equation}
\label{64}
\varphi'' \circ \bar{\iota}(\be_{F_{i}}) \le \varphi'(\be_{F_{i}}) \; \forall \, i =1, \ldots, s\, .
\end{equation}
By (\ref{38}) we have 
\begin{equation}
\label{70}
\|g(\varphi''\bar{\iota}(\be_{F_{i}})) - g(\varphi'(\be_{F_{i}}))\| < \frac{\delta^{2}}{12} \; \forall \, i= 1, \ldots, s \, .
\end{equation}

Moreover, if $\varphi'$ is $n$-decomposable with respect to the decomposition $F=F^{(0)} \oplus \ldots \oplus F^{(n)}$, then we may assume $\varphi''$ to be $n$-decomposable with respect to the decomposition $\bar{F}=\bar{F}^{(0)} \oplus \ldots \oplus \bar{F}^{(n)}$, where $\bar{F}^{(j)}=\bar{\iota}(F^{(j)})$, $j=0, \ldots,n$. In particular, we have 
\begin{equation}
\label{45}
\ord(\varphi''|_{\bar{\iota}(F_{i})}) = 0 \; \forall \, i = 1, \ldots, s \, .
\end{equation}
We denote the matrix blocks of $\bar{F}$ by $\bar{F}_{i}$, $i=1, \ldots, \bar{s}$, and set $\varphi''_{i}:= \varphi''|_{\bar{F}_{i}}$. Each $\varphi''_{i}$ is a multiple of a $*$-homomorphism 
\[
\sigma''_{i}: \bar{F}_{i} \to (\be_{A} - q)A (\be_{A} - q) \, ,
\]
that is, 
\begin{equation}
\label{39}
\varphi''_{i} = \lambda_{i} \cdot \sigma''_{i}
\end{equation}
for some $0 \le \lambda_{i} \le 1$, $i = 1, \ldots, \bar{s}$. Set
\[
\bar{\mu}:= \halb \left(\mu + \frac{1}{2(n+1)} \right) \, .
\]
With $\varrho : \C^{n+1} \to \Zh$ as in Lemma \ref{tracial-division} (using $\bar{\mu}$ in place of $\mu$),  following \cite{W5}, Lemma 2.5, we may define a c.p.\ map
\begin{equation}
\label{60}
\bar{\varphi}:\bar{F} \to (\be_{A}-q)A(\be_{A}-q) \otimes \Zh
\end{equation}
by 
\begin{equation}
\label{40}
\bar{\varphi}(x):= \sum_{j=0}^{n} \varphi''(x \be_{\bar{F}^{(j)}}) \otimes \varrho(e_{j+1}),
\end{equation}
where $e_{1}, \ldots, e_{n+1}$ denote the canonical generators of $\C^{n+1}$. It is obvious that $\bar{\varphi}$ is in fact c.p.c.\ and has order zero, since the $\varrho(e_{j})$ are pairwise orthogonal. By \ref{tracial-division} we have
\begin{equation}
\label{12}
\bar{\tau}(\varrho(e_{j})) > \bar{\mu}
\end{equation}
for $j= 1, \ldots,n+1$, where $\bar{\tau}$ denotes the unique tracial state on $\Zh$. 

For later use we also note that $\bar{\varphi}_{i}:=\bar{\varphi}|_{\bar{F}_{i}}$ satisfies
\begin{eqnarray}
\label{42}
\bar{\varphi}_{i}(x) & \stackrel{(\ref{39},\ref{40})}{=} & \sum_{j=0}^{n} \lambda_{i} \cdot \sigma_{i}''(x \be_{\bar{F}^{(j)}}) \otimes \varrho(e_{j+1}) \nonumber \\
& = & \sigma_{i}''(x) \otimes (\lambda_{i} \cdot \varrho(e_{\bj (i)+1}))  
\end{eqnarray}
for all $x \in \bar{F}_{i}$, $i=1, \ldots, \bar{s}$, where $\bj(i)$ denotes the (uniquely determined) $j \in \{0, \ldots ,n\}$ for which $\be_{\bar{F}_{i}}\be_{\bar{F}^{(j)}} \neq 0$. In particular, we have
\begin{equation}
\label{66}
\bar{\varphi}_{i}(\be_{\bar{F}_{i}}) = \sigma_{i}'' (\be_{\bar{F}_{i}}) \otimes (\lambda_{i} \cdot \varrho(e_{\bj(i) +1})) \, ;
\end{equation}
since $\sigma_{i}''(\be_{\bar{F}_{i}})$ is a projection, it is straightforward to check that
\begin{equation}
\label{2}
f(\bar{\varphi}_{i}(\be_{\bar{F}_{i}})) = \sigma_{i}''(\be_{\bar{F}_{i}}) \otimes f(\lambda_{i} \cdot \varrho(e_{\bj(i) + 1}))
\end{equation}
for any $f \in \Ch_{0}((0,1])$. For $g$ and $h$ defined as above we obtain
\begin{eqnarray}
\sigma_{i}''(x) \otimes g(\lambda_{i} \cdot \varrho(e_{\bj(i)+1})) & \stackrel{(\ref{41})}{=} & \sigma_{i}''(x) \otimes h(\lambda_{i} \cdot \varrho(e_{\bj(i)+1})) (\lambda_{i} \cdot \varrho(e_{\bj(i)+1}))  \nonumber \\
& \stackrel{(\ref{2},\ref{42})}= & h(\bar{\varphi}_{i}(\be_{\bar{F}_{i}})) \bar{\varphi}_{i}(x)  \nonumber \\
& = & \bar{\varphi}_{i}(x) h (\bar{\varphi}_{i}(\be_{\bar{F}_{i}}))  \, . \label{3}
\end{eqnarray}

Choose $\beta'>0$ such that
\begin{equation}
\label{47a}
\bar{s} (4 \beta'(1+\frac{1}{\eta})) + 2 \bar{s}^{2}(\beta')^{\halb}   <   \frac{\delta}{8} \mbox{ and } \beta' < \frac{\delta^{2}\eta^{2}}{32} \, .
\end{equation}
By Proposition \ref{rr0dr} in connection with Proposition \ref{polynomial-appr} (with $\beta'$ in place of $\delta$, $\beta''$ in place of $\beta$ and both $g$ and $h$ in place of $f$) there are a finite-dimensional $C^{*}$-algebra $\tilde{F}$ with a centered embedding $\tilde{\iota}:\bar{F} \to \tilde{F}$ and a c.p.c.\ discrete order zero map
\[
\tilde{\varphi}: \tilde{F} \to (\be_{A} - q) A (\be_{A} - q) \otimes \Zh
\]
such that 
\begin{equation}
\label{9}
\tilde{\varphi}\tilde{\iota}(\be_{\bar{F}_{i}}) \le \bar{\varphi}(\be_{\bar{F}_{i}}) \, ,
\end{equation}
\begin{equation}
\label{7'}
\|\tilde{\varphi}\tilde{\iota}(x) - \bar{\varphi}(x)\| < \beta'' \|x\| < \beta' \|x\| \; \forall \, 0 \neq x \in \bar{F}_{+} \, ,
\end{equation}
\begin{equation}
\label{7}
\|g(\tilde{\varphi}\tilde{\iota}(\be_{\bar{F}_{i}})) - g(\bar{\varphi}(\be_{\bar{F}_{i}}))\|< \beta'
\end{equation}
and
\begin{equation}
\label{8}
\|h(\tilde{\varphi}\tilde{\iota}(\be_{\bar{F}_{i}})) - h(\bar{\varphi}(\be_{\bar{F}_{i}}))\|< \beta'
\end{equation}
for $i=1, \ldots, \bar{s}$.

Let $\chi_{(\eta,1]}$ denote the characteristic function on the interval $(\eta,1]$ and set
\begin{equation}
\label{43}
\bar{p}_{i} := \chi_{(\eta,1]}(\tilde{\varphi} \tilde{\iota}(\be_{\bar{F}_{i}})) \in (\be_{A} - q)A(\be_{A} - q) \otimes \Zh
\end{equation}
for $i=1, \ldots, \bar{s}$; note that the $\bar{p}_{i}$ are well-defined projections in $(\be_{A}-q)A(\be_{A}-q)\otimes \Zh$, since $\tilde{\varphi}$ is a discrete order zero map (whence $\chi_{(\eta,1]}$ is continuous on the spectrum of $\tilde{\varphi}\tilde{\iota}(\be_{\bar{F}_{i}})$ for each $i$). Moreover, the $\bar{p}_{i}$ are pairwise orthogonal (again since $\ord \tilde{\varphi} = 0$), so they add up to a projection 
\begin{equation}
\label{74}
p:= \sum_{i=1}^{\bar{s}} \bar{p}_{i} \, ;
\end{equation}
it is clear that 
\begin{equation}
\label{11}
p = \chi_{(\eta,1]}(\tilde{\varphi}(\be_{\tilde{F}})) \in C^{*}(\tilde{\varphi}(\be_{\tilde{F}})) \subset (\be_{A}-q)A(\be_{A}-q)\otimes \Zh 
\end{equation}
and that
\begin{equation}
\label{56}
p \stackrel{(\ref{41})}{=} p h(\tilde{\varphi} \tilde{\iota} (\be_{\bar{F}})) \tilde{\varphi}\tilde{\iota}(\be_{\bar{F}}) \, . 
\end{equation}
From \cite{W4}, 1.2, we see that $p$ commutes with $\tilde{\varphi}(\tilde{F})$ and that $p \tilde{\varphi}(\, . \, ) = p \tilde{\varphi}( \, . \, )p$ is an order zero map. Define a map $\tilde{\sigma}: \tilde{F} \to (\be_{A}-q)A(\be_{A}-q) \otimes \Zh$ by
\begin{equation}
\label{10}
\tilde{\sigma}( \, . \,):= (p \tilde{\varphi}(\be_{\tilde{F}})p)^{-1} \tilde{\varphi}(\,.\,) \, ,
\end{equation}
where the inverse is well-defined if taken in $p C^{*}(\tilde{\varphi} (\be_{\tilde{F}}))p$. It is obvious  that $\tilde{\sigma}$ is a supporting $*$-homomorphism (in the sense of \cite{W4}, 1.2) for the c.p.c.\ map $p\tilde{\varphi}(\,.\,)p$, i.e., 
\begin{equation}
\label{58}
p \tilde{\varphi}(\, . \,)p = p \tilde{\varphi}(\be_{\tilde{F}}) p \tilde{\sigma}( \, . \,) \, ,
\end{equation}
and that
\begin{equation}
\label{59}
\tilde{\sigma}(\, . \,) = p \tilde{\sigma}(\, . \, ) p \, .
\end{equation}
For $0 \le x \in \Bh_{1}(\bar{F}_{i})$, $i = 1, \ldots, \bar{s}$, we now compute
\begin{eqnarray}
\lefteqn{\|[\bar{p}_{i}, \varphi_{i}''(x) \otimes \be_{\Zh}]\| }\nonumber \\
& \stackrel{(\ref{39})}{=} & |\lambda_{i}| \|[\bar{p}_{i},\sigma_{i}''(x) \otimes \be_{\Zh}]\| \nonumber \\
& \stackrel{(\ref{43},\ref{44})}{=} & |\lambda_{i}| \|\bar{p}_{i} g(\tilde{\varphi} \tilde{\iota} (\be_{\bar{F}_{i}})) (\sigma_{i}''(x) \otimes \be_{\Zh}) - (\sigma_{i}''(x) \otimes \be_{\Zh}) g(\tilde{\varphi} \tilde{\iota} (\be_{\bar{F}_{i}})) \bar{p}_{i}\|  \nonumber \\
& \stackrel{(\ref{7})}{\le} & |\lambda_{i}| \|\bar{p}_{i} g(\bar{\varphi}_{i}(\be_{\bar{F}_{i}}))(\sigma_{i}''(x) \otimes \be_{\Zh}) -  (\sigma_{i}''(x) \otimes \be_{\Zh}) g(\bar{\varphi}_{i}(\be_{\bar{F}_{i}})) \bar{p}_{i} \| +2 \beta'  \nonumber \\
& \stackrel{(\ref{2})}{=} & |\lambda_{i}| \| [\bar{p}_{i}, (\sigma_{i}''(x) \otimes g(\lambda_{i} \cdot \varrho(e_{\bj(i)+1})))] \| + 2 \beta' \nonumber \\
& \stackrel{(\ref{3})}{=} & |\lambda_{i}| \| [\bar{p}_{i}, h(\bar{\varphi}_{i}(\be_{\bar{F}_{i}})) \bar{\varphi}_{i}(x)] \| + 2 \beta' \nonumber \\
& \stackrel{(\ref{8})}{\le} & |\lambda_{i}| \| [\bar{p}_{i}, h(\tilde{\varphi} \tilde{\iota} (\be_{\bar{F}_{i}})) \bar{\varphi}_{i}(x)] \| + 2 \beta' + 2 \beta'  \nonumber \\
& \stackrel{(\ref{7'},\ref{41})}{\le} & |\lambda_{i}| \|[\bar{p}_{i}, h(\tilde{\varphi} \tilde{\iota}(\be_{\bar{F}_{i}})) \tilde{\varphi} \tilde{\iota}(\be_{\bar{F}_{i}}x)]\| + 4 \beta' + \frac{\beta'}{\eta} \nonumber \\
& = & 4 \beta' + \frac{\beta'}{\eta} \, ,  \label{6}
\end{eqnarray}
where for the last equation we have used that $\tilde{\varphi} \tilde{\iota}|_{\bar{F}_{i}}$ is an order zero map, whence the elements of $C^{*}(\tilde{\varphi}\tilde{\iota}(\be_{\bar{F}_{i}}))$ commute with those of $\tilde{\varphi}\tilde{\iota}(\bar{F}_{i})$ for each $i$ (cf.\ \cite{W4}, 1.2).
 
Next note that, for $i = 1, \ldots,\bar{s}$,
\begin{eqnarray*}
\bar{p}_{i} & \stackrel{(\ref{43},\ref{44})}{\le} & g(\tilde{\varphi} \tilde{\iota}(\be_{\bar{F}_{i}})) \\
& \stackrel{(\ref{7},\ref{41})}{\le} & g(\bar{\varphi}_{i} (\be_{\bar{F}_{i}})) + \beta' \cdot \be_{A} \otimes \be_{\Zh} \\
& \stackrel{(\ref{2})}{\le} & \sigma_{i}''(\be_{\bar{F}_{i}}) \otimes \be_{\Zh} + \beta' \cdot \be_{A} \otimes \be_{\Zh} \, .
\end{eqnarray*}
Therefore, if $\varphi_{i}'' \perp \varphi_{j}''$ for some $i, j \in \{1, \ldots , \bar{s}\}$, we have
\begin{eqnarray}
\label{5}
\lefteqn{\| \bar{p}_{i} (\varphi_{j}''(x) \otimes \be_{\Zh})\| } \nonumber \\
& \le & \|(\varphi_{j}''(x) \otimes \be_{\Zh}) \bar{p}_{i} (\varphi_{j}''(x) \otimes \be_{\Zh})\|^{\halb} \nonumber \\
& \stackrel{(\ref{39})}{\le} & \|(\varphi_{j}''(x) \otimes \be_{\Zh}) \sigma_{i}''(\be_{\bar{F}_{i}})(\varphi_{j}''(x) \otimes \be_{\Zh}) + \beta' \cdot (\varphi_{j}''(x) \otimes \be_{\Zh})^{2}\|^{\halb} \nonumber \\
& \le & (\beta')^{\halb} \|x \| \; \forall \, 0 \neq x \in (\bar{F}_{j})_{+} \, .
\end{eqnarray}
For $i=1, \ldots,s$, define
\begin{equation}
\label{63}
I(i) := \{j \in \{1, \ldots, \bar{s} \} \, | \, \be_{\bar{F}_{j}} \le \bar{\iota}(\be_{F_{i}}) \} 
\end{equation}
and 
\begin{equation}
\label{47}
p_{i}' := \sum_{j \in I(i)} \bar{p}_{j} \, ;
\end{equation}
we have
\begin{equation}
\label{76}
\sum_{i=1}^{s} p_{i}' = p \, .
\end{equation}
Note that if $j \neq k \in I(i)$, then
\begin{equation}
\label{46}
\varphi_{j}'' \perp \varphi_{k}'' \, ,
\end{equation}
since $\varphi''|\bar{\iota}(F_{i})$ has order zero for all $i=1, \ldots,s$ by (\ref{45}).

For any $0 \neq x \in (F_{i})_{+}$ and $i=1, \ldots,s$, from (\ref{4}) we obtain
\[
\|\varphi_{i}'(x) - \sum_{j \in I(i)} \varphi_{j}'' \circ \bar{\iota}(x) \|< \beta \|x\| \, ,
\]
whence
\begin{eqnarray}
\label{49}
\lefteqn{\|[p_{i}' , \varphi_{i}' (x) \otimes \be_{\Zh}] \|} \nonumber \\
& < & \| [\sum_{j \in I(i)} \bar{p}_{j}, \sum_{j \in I(i)} \varphi_{j}'' (\be_{\bar{F}_{j}} \bar{\iota}(x)) \otimes \be_{\Zh}] \| + 2 \beta \|x\| \nonumber \\
& \stackrel{(\ref{46},\ref{5})}{\le} & \| \sum_{j \in I(i)} [\bar{p}_{j},  \varphi_{j}'' (\be_{\bar{F}_{j}} \bar{\iota}(x)) \otimes \be_{\Zh}] \| +(2 \bar{s}^{2} (\beta')^{\halb} + 2 \beta) \|x\| \nonumber \\
& \stackrel{(\ref{6})}{\le} &  (\bar{s}(4 \beta' (1+\frac{1}{\eta})) + 2 \bar{s}^{2} (\beta')^{\halb} + 2 \beta) \|x\| \nonumber \\
& \stackrel{(\ref{47a},\ref{48})}{<} & \frac{\delta}{4} \|x\|  \, .
\end{eqnarray} 
Furthermore,
\begin{eqnarray}
\label{17}
\lefteqn{\|[p_{i}' , \varphi_{i}(x) \otimes \be_{\Zh}]\| } \nonumber \\
& \stackrel{(\ref{43},\ref{47})}{=} & \| p_{i}' ((\be_{A}-q) \varphi_{i}(x)) \otimes \be_{\Zh} - (\varphi_{i}(x) (\be_{A} - q)) \otimes \be_{\Zh} p_{i}'\| \nonumber \\
& \stackrel{(\ref{37})}{<} & \|[p_{i}', ((\be_{A}-q) \varphi_{i}(x)(\be_{A}-q)) \otimes \be_{\Zh}] \| + 2 \gamma \|x\| \nonumber \\
& \stackrel{(\ref{13})}{\le} & \|[p_{i}' , \varphi_{i}'(x) \otimes \be_{\Zh}] \| + 2 \gamma \|x\| + 2 \beta \|x\| \nonumber \\
& \stackrel{(\ref{49},\ref{50},\ref{48})}{<} & \left( \frac{\delta}{4} + \frac{\delta}{8} + \frac{\delta}{8} \right) \|x\| \nonumber \\
& = & \frac{\delta}{2} \|x\| \; \forall \, 0 \neq x \in (F_{i})_{+}, \, i=1, \ldots, s \, . 
\end{eqnarray}
Next we check that, for $b \in \bar{\Gh}$ ($=\Gh \cup \{a^{2} \, | \, a \in \Gh\})$,
\begin{eqnarray}
\label{52}
\lefteqn{\|\varphi''\bar{\iota}\psi (b) + \varphi^{\times}\psi(b) - b\| } \nonumber \\
& \stackrel{(\ref{4})}{<} & \|\varphi'\psi(b)+\varphi^{\times}\psi(b) - b \| + \beta \nonumber  \\
& \stackrel{(\ref{13})}{<} & \| (\be_{A}-q)\varphi \psi(b) (\be_{A}-q) + q \varphi \psi(b) q - b\| + 3\beta \nonumber  \\
& \stackrel{(\ref{37})}{<} & \|\varphi \psi(b) - b \| + 3\beta + 2 \gamma \nonumber \\
& \stackrel{(\ref{51})}{<} & \frac{\eta^{6}}{(n+1)^{2}}  + 3\beta + 2 \gamma \nonumber \\
& \stackrel{(\ref{48},\ref{50})}{<} & 2 \frac{\eta^{6}}{(n+1)^{2}} \, .
\end{eqnarray}
From (\ref{52}) and Lemma \ref{multiplicative-domain} (with $(\bar{F} \oplus F, \bar{\iota} \psi \oplus \psi, \varphi'' + \varphi^{\times})$ in place of $(F,\psi,\varphi)$) we see that
\begin{equation}
\label{53}
\|\varphi''(\be_{\bar{F}^{(j)}}) \varphi''\bar{\iota}\psi(b) - \varphi''(\be_{\bar{F}^{(j)}} \bar{\iota}\psi(b))\| < 2 \cdot 2^{\halb} \frac{\eta^{3}}{n+1} \; \forall \, b \in \Gh, \, j=0, \ldots, n \, .
\end{equation}
Since the $\varrho(e_{j})$ are pairwise orthogonal, we even have
\begin{eqnarray}
\label{57}
\lefteqn{\|\bar{\varphi}(\be_{\bar{F}})(\varphi'' \bar{\iota}\psi(b) \otimes \be_{\Zh}) - \bar{\varphi}\bar{\iota} \psi(b)\|} \nonumber \\
& \stackrel{(\ref{40})}{=} & \|\sum_{j=0}^{n} (\varphi''(\be_{\bar{F}^{(j)}}) \varphi''\bar{\iota} \psi(b) - \varphi''(\be_{\bar{F}^{(j)}} \bar{\iota} \psi(b))) \otimes \varrho(e_{j+1}) \| \nonumber \\
& \stackrel{(\ref{53})}{<} & 4 \frac{\eta^{3}}{n+1} \; \forall \, b \in \Gh \, .
\end{eqnarray}
We are now prepared to compute
\begin{eqnarray}
\lefteqn{\|p(b \otimes \be_{\Zh})p - \tilde{\sigma} \tilde{\iota} \bar{\iota} \psi(b) \|} \nonumber \\
& \stackrel{(\ref{51},\ref{54})}{\le} & \|p((\varphi'\psi(b) + \varphi^{\times}\psi(b)) \otimes \be_{\Zh})p - \tilde{\sigma}\tilde{\iota}\bar{\iota}\psi(b)\| + \frac{\eta^6}{(n+1)^2} + \beta \nonumber \\
& \stackrel{(\ref{55},\ref{11})}{=} & \|p (\varphi'\psi(b) \otimes \be_{\Zh})p - \tilde{\sigma}\tilde{\iota}\bar{\iota}\psi(b)\| + \frac{\eta^6}{(n+1)^2} + \beta \nonumber\\
& \stackrel{(\ref{4})}{\le} & \|p(\varphi''\bar{\iota}\psi(b) \otimes \be_{\Zh})p - \tilde{\sigma}\tilde{\iota}\bar{\iota}\psi(b)\| + \frac{\eta^6}{(n+1)^2} + 2 \beta \nonumber \\
& \stackrel{(\ref{56})}{=} & \|p h(\tilde{\varphi}\tilde{\iota}(\be_{\bar{F}})) \tilde{\varphi}\tilde{\iota}(\be_{\bar{F}}) (\varphi'' \bar{\iota}\psi(b) \otimes \be_{\Zh})p - \tilde{\sigma}\tilde{\iota}\bar{\iota}\psi(b)\| + \frac{\eta^6}{(n+1)^2} + 2 \beta  \nonumber \\
& \stackrel{(\ref{7'},\ref{41})}{\le} & \|p h(\bar{\varphi}(\be_{\bar{F}})) \bar{\varphi}(\be_{\bar{F}}) (\varphi'' \bar{\iota}\psi(b) \otimes \be_{\Zh})p - \tilde{\sigma}\tilde{\iota}\bar{\iota}\psi(b)\| \nonumber \\
& &  + \frac{\eta^6}{(n+1)^2} + 2 \beta  + \beta' + 2 \frac{\beta'}{\eta} \nonumber \\
& \stackrel{(\ref{57},\ref{41})}{\le} & \|p h(\bar{\varphi}(\be_{\bar{F}})) \bar{\varphi} \bar{\iota}\psi(b)p   - \tilde{\sigma}\tilde{\iota}\bar{\iota}\psi(b)\| \nonumber \\
& & + \frac{\eta^6}{(n+1)^2} + 2 \beta  + \beta' + 2 \frac{\beta'}{\eta} + 12 \frac{\eta^{2}}{n+1} \nonumber \\
& \stackrel{(\ref{7'},\ref{41})}{\le} & \|p h(\tilde{\varphi}\tilde{\iota}(\be_{\bar{F}})) \tilde{\varphi}\tilde{\iota} \bar{\iota}\psi(b)p   - \tilde{\sigma}\tilde{\iota}\bar{\iota}\psi(b)\| \nonumber \\
& & + \frac{\eta^6}{(n+1)^2} + 2 \beta  + \beta' + 2 \frac{\beta'}{\eta} + 12 \frac{\eta^{2}}{n+1} + \beta' + 2 \frac{\beta'}{\eta} \nonumber \\
& \stackrel{(\ref{58})}{=} & \|p h(\tilde{\varphi}\tilde{\iota}(\be_{\bar{F}})) \tilde{\varphi}\tilde{\iota} (\be_{\bar{F}}) \tilde{\sigma}\tilde{\iota} \bar{\iota}\psi(b)p   - \tilde{\sigma}\tilde{\iota}\bar{\iota}\psi(b)\| \nonumber \\
& & + \frac{\eta^6}{(n+1)^2} + 2 \beta  + \beta' + 2 \frac{\beta'}{\eta} + 12 \frac{\eta^{2}}{n+1} + \beta' + 2 \frac{\beta'}{\eta} \nonumber \\
& \stackrel{(\ref{56},\ref{59})}{=} & 0+\frac{\eta^6}{(n+1)^2} + 2 \beta + 2 \beta' + 4 \frac{\beta'}{\eta} + 12 \frac{\eta^{2}}{n+1} \nonumber \\
& \stackrel{(\ref{48},\ref{47a})}{<} & \frac{3}{4} \eta \label{15}
\end{eqnarray}
for all $b \in \Gh$.

If $\tau \in T(A)$ is a tracial state, then
\begin{eqnarray}
\lefteqn{\tau \otimes \bar{\tau} (p)} \nonumber \\
& \stackrel{(\ref{11})}{\ge} & \tau \otimes \bar{\tau} (\tilde{\varphi}(\be_{\tilde{F}})) - \eta \cdot \tau(\be_{A}-q) \nonumber \\
& \stackrel{(\ref{7'},\ref{60})}{\ge} & \tau \otimes \bar{\tau} (\bar{\varphi}(\be_{\bar{F}})) - (\eta + \beta') \cdot \tau(\be_{A}-q) \nonumber \\
& \stackrel{(\ref{42})}{=} & \sum_{j=0}^{n} \tau(\varphi''(\be_{\bar{F}^{(j)}})) \bar{\tau}(\varrho(e_{\bj(i)+1}))     - (\eta + \beta') \cdot \tau(\be_{A}-q) \nonumber \\
& \stackrel{(\ref{12})}{\ge} & \bar{\mu} \cdot \sum_{j=0}^{n} \tau(\varphi''(\be_{\bar{F}^{(j)}})) - (\eta + \beta') \cdot \tau(\be_{A}-q) \nonumber \\
& \stackrel{(\ref{4},\ref{61})}{\ge} & \bar{\mu} \cdot \tau (\varphi'(\be_{F})) -  (\beta + \eta + \beta') \cdot \tau(\be_{A}-q) \nonumber \\
& \stackrel{(\ref{13})}{\ge} & \bar{\mu} \cdot \tau((\be_{A}-q) \varphi(\be_{F})(\be_{A}-q)) - (2 \beta + \eta + \beta') \cdot \tau(\be_{A}-q) \nonumber \\
& \stackrel{(\ref{51})}{\ge} & \bar{\mu} \cdot \tau(\be_{A}-q) - (\eta + 2 \beta + \eta + \beta') \cdot \tau(\be_{A}-q) \nonumber \\
& \stackrel{(\ref{48},\ref{47a})}{\ge} & (\bar{\mu} - 4 \eta) \cdot \tau(\be_{A}-q) \nonumber \\
& \stackrel{(\ref{62})}{>} & (\mu + \eta) \cdot \tau(\be_{A}-q) \, . \label{16}
\end{eqnarray}

For $i \in \{1, \ldots,s\}$ we have 
\begin{eqnarray}
\label{65}
\sum_{j \in I(i)} \tilde{\varphi}\tilde{\iota} (\be_{\bar{F}_{j}}) & \stackrel{(\ref{9})}{\le} & \sum_{j \in I(i)} \bar{\varphi}(\be_{\bar{F}_{j}}) \nonumber \\
& \stackrel{(\ref{40})}{\le} & \sum_{j \in I(i)} (\sum_{k=0}^{n} \varphi''(\be_{\bar{F}_{j}} \be_{\bar{F}^{(k)}}) \otimes \be_{\Zh}) \nonumber \\
& \stackrel{(\ref{63})}{=} & \varphi''(\bar{\iota}(\be_{F_{i}})) \otimes \be_{\Zh} \nonumber \\
& \stackrel{(\ref{64})}{\le} & \varphi'(\be_{F_{i}}) \otimes \be_{\Zh} \, ,
\end{eqnarray}
whence 
\begin{eqnarray*}
p_{i}' & \stackrel{(\ref{47})}{=} & \sum_{j \in I(i)} \bar{p}_{j} \\
& \stackrel{(\ref{43})}{=} & \sum_{j \in I(i)} \chi_{(\eta,1]}(\tilde{\varphi}\tilde{\iota}(\be_{\bar{F}_{j}})) \in \overline{(\varphi'(\be_{F_{i}})\otimes \be_{\Zh})(A \otimes \Zh) (\varphi'(\be_{F_{i}}) \otimes \be_{\Zh})} \, .
\end{eqnarray*}
Even more, 
\begin{eqnarray*}
\eta \cdot p_{i}' & \le &  \sum_{j \in I(i)} \tilde{\varphi}\tilde{\iota}(\be_{\bar{F}_{i}}) \\
& \stackrel{(\ref{65})}{\le} &  \varphi'(\be_{F_{i}}) \otimes \be_{\Zh} \\
& \stackrel{(\ref{13})}{\le} &  ((\be_{A}-q)\varphi(\be_{F_{i}})(\be_{A}-q) \otimes \be_{\Zh} + q \varphi(\be_{F_{i}})q \otimes \be_{\Zh}) + \beta \cdot \be_{A} \otimes \be_{\Zh} \\
& \stackrel{(\ref{37})}{\le} & (2 \gamma + \beta) \cdot \be_{A} \otimes \be_{\Zh} +  \varphi(\be_{F_{i}}) \otimes \be_{\Zh} \, ,
\end{eqnarray*}
and it follows from Proposition \ref{almost-hereditary} that
\begin{equation}
\label{14}
\dist(p_{i}', \overline{(\varphi(\be_{F_{i}}) \otimes \be_{\Zh})A \otimes \Zh(\varphi(\be_{F_{i}}) \otimes \be_{\Zh})}) \le \frac{3}{\eta} (2\gamma + \beta)^{\halb} \stackrel{(\ref{48},\ref{50})}{<} \frac{\delta}{2s} \, .
\end{equation}

From (\ref{7}) and the fact that $\tilde{\varphi}\tilde{\iota}$ is subordinate (in the sense of \cite{W4}, Definition 1.4) to the order zero map $\bar{\varphi}$ we know that, for $i \in \{1, \ldots, s \}$,
\begin{eqnarray}
\label{68}
\|g(\tilde{\varphi}\tilde{\iota}\bar{\iota}(\be_{F_{i}})) - g(\bar{\varphi}\bar{\iota}(\be_{F_{i}}))\| & \le & \max_{j \in I(i)} \|g(\tilde{\varphi}\tilde{\iota}(\be_{\bar{F}_{j}})) - g(\tilde{\varphi}(\be_{\bar{F}_{j}}))\| \nonumber \\
& \stackrel{(\ref{32})}{<} & \beta' \, .
\end{eqnarray}
For $i \in \{1, \ldots, \bar{s}\}$ we have 
\begin{eqnarray}
\label{67}
g(\bar{\varphi}(\be_{\bar{F}_{i}})) & \stackrel{(\ref{66})}{=} & g(\sigma_{i}''(\be_{\bar{F}_{i}}) \otimes \lambda_{i} \cdot \varrho(e_{\bj(i)+1})) \nonumber \\
& \stackrel{(\ref{2})}{=} & \sigma_{i}'' (\be_{\bar{F}_{i}}) \otimes g(\lambda_{i} \cdot \varrho(e_{\bj(i)+1})) \nonumber  \\
& \le & \sigma_{i}''(\be_{\bar{F}_{i}}) \otimes g(\lambda_{i}) \cdot \be_{\Zh} \nonumber \\
& = & g(\sigma_{i}''(\be_{\bar{F}_{i}}) \otimes \lambda_{i} \cdot \be_{\Zh}) \nonumber \\
& \stackrel{(\ref{2})}{=} & g(\varphi_{i}''(\be_{\bar{F}_{i}})) \otimes \be_{\Zh} \, ,
\end{eqnarray}
where the inequality follows from the fact that $g(\lambda \cdot t) \le g(\lambda)$ for all $0\le \lambda,t \le 1$: the latter implies that the constant function $g(\lambda) \cdot \be_{[0,1]}$ on $[0,1]$ dominates the function $(t \mapsto g(\lambda \cdot t)) \in \Ch([0,1])$; Gelfand's theorem now yields $g(\lambda \cdot a) \le g(\lambda) \cdot \be$ for any $0\le \lambda \le 1$ and any $0\le a \le \be$ in a unital $C^{*}$-algebra. \\
Since $\bar{\varphi}$ (by construction) and  $\varphi''|_{\bar{\iota}(F_{i})}$ (by (\ref{45})) have order zero, we even have 
\begin{eqnarray}
\label{69}
g(\bar{\varphi}\bar{\iota}(\be_{F_{i}})) & = & \sum_{j \in I(i)} g(\bar{\varphi}(\be_{\bar{F}_{j}})) \nonumber \\
& \stackrel{(\ref{67})}{\le} & \sum_{j \in I(i)} g(\varphi''_{i}(\be_{\bar{F}_{i}})) \otimes \be_{\Zh} \nonumber \\
& = &  g(\varphi''\bar{\iota}(\be_{F_{i}})) \otimes \be_{\Zh} \; \forall \, i \in \{1, \ldots, s\} \, .
\end{eqnarray}
We conclude that
\begin{eqnarray}
\label{72}
g(\tilde{\varphi}\tilde{\iota} \bar{\iota}(\be_{F_{i}})) & \stackrel{(\ref{68})}{\le} & g(\bar{\varphi}\bar{\iota}(\be_{F_{i}})) + \beta' \cdot \be_{A} \otimes \be_{\Zh} \nonumber \\
& \stackrel{(\ref{69})}{\le} & g(\varphi''\bar{\iota}(\be_{F_{i}})) \otimes \be_{\Zh} + \beta' \cdot \be_{A} \otimes \be_{\Zh} \nonumber \\
& \stackrel{(\ref{70})}{\le} & g(\varphi'(\be_{F_{i}})) \otimes \be_{\Zh} + (\beta' + \frac{\delta^{2}}{12}) \cdot \be_{A}\otimes \be_{\Zh} \nonumber \\
& \stackrel{(\ref{71})}{\le} & g(\varphi(\be_{F_{i}})) \otimes \be_{\Zh} + (\beta' + \frac{\delta^{2}}{12} + \beta) \cdot \be_{A}\otimes \be_{\Zh} \, .
\end{eqnarray} 
for $i \in \{1, \ldots, s\}$. Now since 
\[
p_{i}' \stackrel{(\ref{47},\ref{43})}{=} g(\tilde{\varphi}\tilde{\iota}\bar{\iota}(\be_{F_{i}}))p_{i}'
\]
for $i \in \{1, \ldots, s\}$, we have
\begin{eqnarray*}
\|p_{i}' (g(\varphi(\be_{F_{i}})) \otimes \be_{\Zh}) - p_{i}'\|^{2} & \le & \|p_{i}' (\be_{A \otimes \Zh} - g(\varphi(\be_{F_{i}})) \otimes \be_{\Zh}) p_{i}' \| \\
& \stackrel{(\ref{72})}{\le} & \beta' + \frac{\delta^{2}}{12} + \beta \\
& \stackrel{(\ref{47a})}{<} & \frac{\delta^{2}}{4} \, ,
\end{eqnarray*}
hence 
\begin{equation}
\label{73}
\|p_{i}' (g(\varphi(\be_{F_{i}})) \otimes \be_{\Zh}) - p_{i}'\| < \frac{\delta}{2} 
\end{equation}
for each $i \in \{1, \ldots , s\}$.

By \cite{TW1}, Remark 2.7, there is a unital $*$-homomorphism $\theta: A \otimes \Zh \to A$ satisfying
\begin{equation}
\label{75}
\|\theta(b \otimes \be_{\Zh}) - b\| < \frac{\eta}{4}  \; \forall \, b \in \Gh \, ,
\end{equation}
\begin{equation}
\label{19}
\| \theta (\varphi(x) \otimes \be_{\Zh}) - \varphi(x)\| < \frac{\delta}{4} \|x\| \; \forall \, 0 \neq x \in F_{+} \, 
\end{equation}
\begin{equation}
\label{20}
\|\theta((\be_{A}-q) \otimes \be_{\Zh}) - (\be_{A}-q)\| < \frac{\eta}{\mu + \eta} \cdot \min_{\tau \in T(A)} \{\tau(\be_{A}-q)\} 
\end{equation}
and
\begin{equation}
\label{21}
\|\theta(g(\varphi(\be_{F_{i}})) \otimes \be_{\Zh}) - g(\varphi(\be_{F_{i}})) \| < \frac{\delta}{2} \; \forall \, i \in \{1, \ldots, s\} 
\end{equation}
(note that $\min_{\tau \in T(A)} \{\tau(\be_{A} - q) \}$ exists and is nonzero since $A$ is unital and simple, whence $T(A)$ is compact and $\tau(\be_{A} - q)>0 \; \forall \, \tau \in T(A)$).\\
Using (\ref{14}), it is straightforward to check that we may even assume that
\begin{equation}
\label{18}
\dist(\theta(p_{i}'),\overline{\varphi(\be_{F_{i}})A\varphi(\be_{F_{i}})}) < \frac{\delta}{s} \; \forall \, i \in \{1, \ldots, s\} \, .
\end{equation}

Define a finite-dimensional $C^{*}$-algebra $C \subset A$ by
\[
C:= \theta \tilde{\sigma}(\tilde{F})
\]
and projections $p_{1}, \ldots, p_{s} \in A$ by 
\begin{equation}
\label{77}
p_{i}:= \theta(p_{i}'), i=1, \ldots,s \, .
\end{equation}
It is clear from our construction that $p_{i} \in C \; \forall \, i$ and that
\[
\sum_{i=1}^{s} p_{i} \stackrel{(\ref{76})}{=} \theta(p) \stackrel{(\ref{59})}{=} \be_{C} \, .
\]
We proceed to check assertions (i), (ii) and (iii) of the lemma:
\begin{eqnarray*}
\lefteqn{\dist(\be_{C}b\be_{C},C)}\\
& \stackrel{(\ref{75})}{\le} & \dist(\theta(p)\theta(b \otimes \be_{\Zh}) \theta(p), \theta(\tilde{\sigma}(\tilde{F}))) + \frac{\eta}{4}  \\
& \le & \|p (b \otimes \be_{\Zh}) p - \tilde{\sigma} \tilde{\iota} \bar{\iota} \psi(b) \| +  \frac{\eta}{4}  \\
& \stackrel{(\ref{15})}{<} & \eta \; \forall \, b \in \Gh \, .
\end{eqnarray*}
If $\tau \in T(A)$, then $\tau \circ \theta \in T(A \otimes \Zh)$, whence there is $\tau' \in T(A)$ such that $\tau \circ \theta = \tau' \otimes \bar{\tau}$. Therefore,
\begin{eqnarray*}
\tau(\be_{C}) & \stackrel{(\ref{76})}{=} & \tau \circ \theta(p) \\
& = & (\tau' \otimes \bar{\tau})(p)\\
& \stackrel{(\ref{16})}{\ge} & (\mu + \eta) \tau'(\be_{A} - q) \\
& = & (\mu + \eta) (\tau' \otimes \bar{\tau})((\be_{A}-q) \otimes \be_{\Zh}) \\
& = & (\mu + \eta) \tau \circ \theta ((\be_{A}-q) \otimes \be_{\Zh}) \\
& \stackrel{(\ref{20})}{>} & (\mu + \eta)(\tau(\be_{A}-q) - \frac{\eta}{\mu + \eta} \tau(\be_{A}-q)) \\
& = & \mu \tau(\be_{A}-q) \, .
\end{eqnarray*}
We also have
\begin{eqnarray*}
\|[p_{i}, \varphi(\be_{F_{i}} x)]\| & \stackrel{(\ref{77},\ref{19})}{\le} & \|[\theta(p_{i}'), \theta(\varphi(\be_{F_{i}} x) \otimes \be_{\Zh})]\| + \frac{\delta}{2} \|x\| \\
& \le & \|[p_{i}', \varphi(\be_{F_{i}} x) \otimes \be_{\Zh}]\| + \frac{\delta}{2} \|x\| \\
& \stackrel{(\ref{17})}{<} & \delta \|x\| 
\end{eqnarray*}
for all $0 \neq x \in F_{+}$,
\[
\dist(p_{i}, \overline{\varphi(\be_{F_{i}})A \varphi(\be_{F_{i}})}) \stackrel{(\ref{18})}{<} \frac{\delta}{s}
\]
and
\begin{eqnarray}
\|p_{i} g(\varphi(\be_{F_{i}})) - p_{i}\| & \stackrel{(\ref{21})}{<} & \|\theta(p_{i}') \theta(g(\varphi(\be_{F_{i}}))) - \theta(p_{i}')\| + \frac{\delta}{2} \nonumber \\
& \le & \|p_{i}' g (\varphi(\be_{F_{i}})) - p_{i}' \| + \frac{\delta}{2} \nonumber \\
& \stackrel{(\ref{73})}{<} & \delta
\end{eqnarray}
for $i=1, \ldots,s$. We are done.    
\end{nproof}
\en

\section{The proof of Theorem \ref{lfdrtr0}}

This section is entirely devoted to the proof of Theorem \ref{lfdrtr0}, following the outline of Section 3. Let $A$ be separable, simple, unital and $\Zh$-stable with real rank zero and locally finite decomposition rank.

Since $A$ has real rank zero, every nonzero hereditary subalgebra contains a nontrivial projection; since $A$ is nuclear and $\Zh$-stable, it satisfies Blackadar's second fundamental comparability property by \cite{R2}, Corollary 4.6. Therefore, it will suffice to show that $A$ satisfies the hypotheses of Proposition \ref{wu-tr0}. 

So let $\varepsilon>0$ and a finite subset $\Fh \subset A$ be given. Without loss of generality we may assume that $\be_{A} \in \Fh$ and that the elements of $\Fh$ are positive and normalized. Moreover, since $A$ has locally finite decomposition rank, we can assume that $\Fh \subset \Bh_{1}(B)_{+}$, where $B \subset A$ is a unital $C^{*}$-subalgebra with $\dr B=n$ for some $n \in \N$. \\
Fix some $0 < \mu < \frac{1}{2(n+1)}$. For $k \in \N$, define 
\begin{equation}
\label{80}
\zeta_{k}:=\mu \sum_{l=0}^{k} (1-\mu)^{l} \, ,
\end{equation}
then
\[
\zeta_{k} \stackrel{k \to \infty}{\longrightarrow} \mu \sum_{l=0}^{\infty} (1-\mu)^{l} = \mu \frac{1}{1-(1-\mu)} = 1 \, ,
\]
whence there is $K \in \N$ such that 
\begin{equation} 
\label{81}
\zeta_{K}>1-\varepsilon \, .
\end{equation}
Define $\Gh_{0}:=\Fh$ and choose $\eta_{0}>0$ such that
\[
\eta_{0}< \min \left\{\frac{\varepsilon}{8}, \, \frac{1}{10}\left( \frac{1}{2(n+1)} - \mu \right) , \, \frac{1}{48}\right\} \, .
\]
Apply Corollary \ref{excisible-appr} (with $\Gh_{0}$ in place of $\Gh$ and $\eta_{0}$ in place of $\eta$) to obtain an $n$-decomposable c.p.\ approximation $(F_{0},\psi_{0},\varphi_{0})$ and $0 < \delta_{0}< \halb$ such that a) and b) of Corollary \ref{excisible-appr} hold. 

Now the hypotheses of Lemma \ref{excision} are fulfilled (with $(F_{0},\psi_{0},\varphi_{0})$, $\eta_{0}$, $\Gh_{0}$ and $\delta_{0}$ in place of $(F,\psi,\varphi)$, $\eta$, $\Gh$ and $\delta$); note that (\ref{51}) is satisfied by Corollary \ref{excisible-appr}a). We obtain $\gamma_{0}>0$ such that the assertion of Lemma \ref{excision} holds.

Next, suppose $\Gh_{k}$, $\eta_{k}$, $(F_{k},\psi_{k},\varphi_{k})$, $\delta_{k}$ and $\gamma_{k}$ have been constructed for some $k \in \N$. Define $\Gh_{k+1}:= \Gh_{k} \cup \varphi_{k}(\Bh_{1}(F_{k})_{+})$ and choose $\eta_{k+1}>0$ such that
\[
\eta_{k+1}<\frac{1}{2^{k+1}} \min\left\{\frac{\varepsilon}{8}, \, \gamma_{k}, \delta_{k}\right\}
\]
and
\[
\eta_{k+1} < \min \left\{\frac{\varepsilon}{8}, \, \frac{1}{10} \left( \frac{1}{2(n+1)} - \mu \right) , \, \frac{1}{48} \right\} \, .
\]
From Corollary \ref{excisible-appr} (with $\Gh_{k+1}$ in place of $\Gh$ and $\eta_{k+1}$ in place of $\eta$) we obtain an $n$-decomposable c.p.\ approximation $(F_{k+1},\psi_{k+1},\varphi_{k+1})$ of $B$ and $0 < \delta_{k+1}< \halb $ such that a) and b) of \ref{excisible-appr} hold. \\
Again, the hypotheses of Lemma \ref{excision} are fulfilled (with $(F_{k+1},\psi_{k+1},\varphi_{k+1})$, $\eta_{k+1}$, $\Gh_{k+1}$ and $\delta_{k+1}$ in place of $(F,\psi,\varphi)$, $\eta$, $\Gh$ and $\delta$), so we obtain $\gamma_{k+1}>0$ such that the assertion of Lemma \ref{excision} holds; we may asume that $\gamma_{k+1}<\gamma_{k}$. 

Induction yields compact subsets $\Gh_{k} \subset B$, positive numbers $\eta_{k}$, $\delta_{k}$, $\gamma_{k}$ and $n$-decomposable c.p.\ approximations $(F_{k},\psi_{k},\varphi_{k})$ for each $k \in \N$. By construction, we have in particular that
\begin{equation}
\label{24}
\sum_{l=0}^{\infty} \eta_{l}< \frac{\varepsilon}{2} \, , \; \sum_{l=k+1}^{K} \eta_{l} < \gamma_{k}
\end{equation}
and 
\[
\Gh_{k} = \Fh \cup \bigcup_{l=0}^{k-1} \varphi_{l}(\Bh_{1}(F_{l})_{+}) \subset \Gh_{k+1} \, .
\]
For each $k$, we denote the summands of $F_{k}$ by $F_{k,i}$, $i=1, \ldots,s_{k}$, in other words, we write $F_{k}= \bigoplus_{i=1}^{k}F_{k,i}$ with matrix algebras $F_{k,i}$. \\
Let $q_{K} \in A$ be the zero projection, then
\[
\|[q_{K},\varphi_{K}(x)]\|=0 < \gamma_{K} \|x\| \; \forall \, 0 \neq x \in (F_{K})_{+}
\]
and by Lemma \ref{excision} there is a finite-dimensional $C^{*}$-subalgebra 
\[
C_{K} \subset (\be_{A}-q_{K}) A (\be_{A}-q_{K}) = A
\]
satisfying
\begin{itemize}
\item[(i)] $\dist(\be_{C_{K}},b,\be_{C_{K}},C_{K})< \eta_{K} \; \forall \, b \in \Gh_{K}$
\item[(ii)] $\tau(\be_{C_{K}}) \ge \mu \cdot \tau(\be_{A}-q_{K})=\mu \; \forall \, \tau \in T(A)$ 
\item[(iii)] the projection $\be_{C_{K}}$ can be written as a sum of $s_{K}$ pairwise orthogonal projections $p_{K,1}, \ldots,p_{K,s_{K}} \in C_{K}$, $\be_{C_{K}}= \sum_{i=1}^{s_{K}}p_{K,i}$, satisfying
\[
\|[p_{K,i}, \varphi_{K}(\be_{F_{K,i}}x)]\| < \delta_{K}\|x\| \; \forall \, 0 \neq x \in (F_{K})_{+} \, ,  
\]
\[
\|p_{K,i} g_{\frac{\eta_{K}}{2},\eta_{K}}(\varphi_{K}(\be_{F_{K,i}})) - p_{K,i}\|<\delta_{K}
\]
and
\[
\dist(p_{K,i} , \overline{\varphi_{K}(\be_{F_{K,i}})A\varphi_{K}(\be_{F_{K,i}})} ) < \frac{\delta_{K}}{s_{K}}
\]
for $i=1, \ldots,s_{K}$.
\end{itemize}
Suppose that, for some $k \in \{1, \ldots,K\}$, we have already constructed pairwise orthogonal finite-dimensional $C^{*}$-subalgebras $C_{l} \subset A$ and projections 
\begin{equation}
\label{79}
q_{l}=\sum_{m=l+1}^{K} \be_{C_{m}} \in A
\end{equation}
for $l=k, \ldots,K$, which satisfy 
\[
C_{l} \subset (\be_{A}-q_{l})A(\be_{A}-q_{l}) \, , \; \|[q_{l},\varphi_{l}(x)]\|< \gamma_{l} \|x\| \; \forall \, 0 \neq x \in F_{l}
\]
and
\begin{itemize}
\item[(i')] $\dist(\be_{C_{l}},b,\be_{C_{l}},C_{l})<  \eta_{l} \; \forall \, b \in \Gh_{l}$
\item[(ii')] $\tau(\be_{C_{l}}) \ge \mu \cdot \tau(\be_{A}-q_{l}) \; \forall \, \tau \in T(A)$ 
\item[(iii')] the projection $\be_{C_{l}}$ can be written as a sum of $s_{l}$ pairwise orthogonal projections $p_{l,1}, \ldots,p_{l,s_{l}} \in C_{l}$, $\be_{C_{l}}= \sum_{i=1}^{s_{l}}p_{l,i}$, satisfying
\[
\|[p_{l,i}, \varphi_{l}(\be_{F_{l,i}}x)]\| < \delta_{l}\|x\| \; \forall \, 0 \neq x \in (F_{l})_{+} \, , 
\]
\[
\|p_{l,i} g_{\frac{\eta_{l}}{2},\eta_{l}}(\varphi_{l}(\be_{F_{l,i}})) - p_{l,i}\|<\delta_{l}
\]
and
\[
\dist(p_{l,i} , \overline{\varphi_{l}(\be_{F_{l,i}})A\varphi_{l}(\be_{F_{l,i}})}) < \frac{\delta_{l}}{s_{l}}
\]
for $i=1, \ldots,s_{l}$.
\end{itemize}
Now (iii') and Corollary \ref{excisible-appr}b) imply that 
\begin{equation}
\label{23}
\|[\be_{C_{l}},b]\| < \eta_{l} \; \forall \, b \in \Gh_{l}, \, l=k, \ldots, K\, .
\end{equation}
Set
\[
q_{k-1}:=q_{k}+\be_{C_{k}} = \sum_{l=k}^{K} \be_{C_{l}} \, ,
\]
then $q_{k-1}$ is a projection since $q_{k} \perp \be_{C_{k}}$ and 
\begin{eqnarray*}
\frac{1}{\|x\|} \cdot \|[q_{k-1},\varphi_{k-1}(x)]\| & \le & \sum_{l=k}^{K} \|[\be_{C_{l}},\frac{1}{\|x\|} \cdot \varphi_{k-1}(x)]\| \\
& \stackrel{(\ref{23})}{<} & \sum_{l=k}^{K} \eta_{l} \\
& \stackrel{(\ref{24})}{\le} & \gamma_{k-1}  \; \forall \, 0\neq x \in (F_{k-1})_{+} \, ,
\end{eqnarray*}
since
\[
\varphi_{k-1}(\Bh_{1}((F_{k-1})_{+})) \subset \Gh_{l} \; \forall \, l=k, \ldots,K\, .
\]
Now by Lemma \ref{excision} there is a finite-dimensional $C^{*}$-subalgebra 
\[
C_{k-1} \subset (\be_{A}-q_{k-1})A(\be_{A}-q_{k-1})
\]
such that
\begin{itemize}
\item[(i'')] $\dist(\be_{C_{k-1}},b,\be_{C_{k-1}},C_{k-1})<  \eta_{k-1} \; \forall \, b \in \Gh_{k-1}$
\item[(ii'')] $\tau(\be_{C_{k-1}}) \ge \mu \cdot \tau(\be_{A}-q_{k-1}) \; \forall \, \tau \in T(A)$ 
\item[(iii'')] the projection $\be_{C_{k-1}}$ can be written as a sum of $s_{k-1}$ pairwise orthogonal projections $p_{k-1,1}, \ldots,p_{k-1,s_{k-1}} \in C_{k-1}$, $\be_{C_{k-1}}= \sum_{i=1}^{s_{k-1}}p_{k-1,i}$, satisfying
\[
\|[p_{k-1,i}, \varphi_{k-1}(\be_{F_{k-1,i}}x)]\| < \delta_{k-1}\|x\| \; \forall \, 0 \neq x \in (F_{k-1})_{+} \, , 
\]
\[
\|p_{k-1,i} g_{\frac{\eta_{k-1}}{2},\eta_{k-1}}(\varphi_{k-1}(\be_{F_{k-1,i}})) - p_{k-1,i}\|<\delta_{k-1}
\]
and
\[
\dist(p_{k-1,i} , \overline{\varphi_{k-1}(\be_{F_{k-1,i}})A\varphi_{k-1}(\be_{F_{k-1,i}})})< \frac{\delta_{k-1}}{s_{k-1}}
\]
for $i=1, \ldots,s_{k-1}$.
\end{itemize}
Induction yields pairwise orthogonal finite-dimensional $C^{*}$-subalgebras $C_{k}\subset A$ and projections $q_{k} \in A$ satisfying $q_{k}= \sum_{m=k+1}^{K} \be_{C_{m}}$ and (i'), (ii') and (iii') above for $k=0, \ldots, K$ in place of $l$. Note that (iii') and Corollary \ref{excisible-appr} b) imply that (\ref{23}) holds for all $l=0,\ldots,K$.\\
Define a finite-dimensional $C^{*}$-subalgebra $D$ of $A$ by 
\[
D:= \bigoplus_{k=0}^{K} C_{k} \, .
\]
We proceed to check properties (i), (ii) and (iii) of Proposition \ref{wu-tr0}. 

First, we have for any $b \in \Fh$
\begin{eqnarray*}
\|[\be_{D},b]\| & \le & \sum_{k=0}^{K} \|[\be_{C_{k}},b]\| \\
& \stackrel{(\ref{23})}{<} & \sum_{k=0}^{K} \eta_{k} \\
& \stackrel{(24)}{<} & \varepsilon \, ,
\end{eqnarray*}
since $\Fh \subset \Gh_{k}$ for $k=0, \ldots, K$ and since (\ref{23}) holds for $l = 0, \ldots, K$. Similarly, we obtain
\begin{eqnarray*}
\lefteqn{\dist(\be_{D}b\be_{D},D)}\\
 & = & \dist((\sum_{k=0}^{K}\be_{C_{k}})(\sum_{k=0}^{K} b \be_{C_{k}}),D) \\
 & \le & \dist((\sum_{k=0}^{K}\be_{C_{k}})(\sum_{k=0}^{K} \be_{C_{k}} b \be_{C_{k}}),D) + \sum_{k=0}^{K} \|[b,\be_{C_{k}}]\| \\
 & \stackrel{(\ref{23})}{<} & \dist(\sum_{k=0}^{K}(\be_{C_{k}}b\be_{C_{k}}),D) + \sum_{k=0}^{K} \eta_{k} \\
 & = & \max_{k=0, \ldots,K} (\dist(\be_{C_{k}}b \be_{C_{k}},C_{k})) + \sum_{k=0}^{K} \eta_{k} \\
 & \stackrel{({\small \mbox{i'}})}{<} & \max_{k=0, \ldots,K}(\eta_{k}) + \sum_{k=0}^{K} \eta_{k} \\
 & \stackrel{(24)}{<} & \varepsilon
\end{eqnarray*}
for any $b \in \Fh$. 

Finally, we show by induction that
\begin{equation}
\label{78}
\tau \left(\sum_{l=0}^{k} \be_{C_{K-l}} \right) \ge \zeta_{k}
\end{equation}
for $k=0, \ldots, K$ and any $\tau \in T(A)$. From Lemma \ref{excision}(ii) we see that
\[
\tau(\be_{C_{K}}) \ge \mu \cdot \tau(\be_{A}-q_{0}) = \mu \cdot \tau(\be_{A}) = \mu = \zeta_{0} \; \forall \, \tau \in T(A) \, ,
\]
so (\ref{78}) holds for $k=0$. Next, suppose we have shown (\ref{78}) for some $k \in \{0, \ldots, K-1\}$ and all $\tau \in T(A)$. Then, 
\begin{eqnarray*}
\lefteqn{\tau \left(\sum_{l=0}^{k+1} \be_{C_{K-l}} \right)} \\
& = & \tau(\be_{C_{K-(k+1)}}) + \tau\left( \sum_{l=0}^{k} \be_{C_{K-l}}\right) \\
& \stackrel{({\small \mbox{ii'}})}{\ge} & \mu \cdot \tau(\be_{A} -q_{K-(k+1)}) +  \tau\left( \sum_{l=0}^{k} \be_{C_{K-l}}\right) \\
& \stackrel{(\ref{79})}{=} & \mu \cdot \tau \left(\be_{A} - \sum_{l=K-(k+1)+1}^{K} \be_{C_{l}} \right)  + \tau\left( \sum_{l=0}^{k} \be_{C_{K-l}}\right) \\
& = & \mu \cdot \tau(\be_{A}) + (1 - \mu) \cdot \tau\left( \sum_{l=0}^{k} \be_{C_{K-l}}\right) \\
& \stackrel{(\ref{78})}{\ge} & \mu + (1-\mu) \zeta_{k}\\
& \stackrel{(\ref{80})}{=} & \zeta_{k+1}
\end{eqnarray*}
for all $\tau \in T(A)$. Therefore, (\ref{78}) holds for all $k=0, \ldots, K$ and $\tau \in T(A)$. In particular,
\[
\tau(\be_{D}) = \tau \left( \sum_{l=0}^{K} \be_{C_{K-l}} \right) \ge \zeta_{K} \stackrel{(\ref{81})}{>} 1 - \varepsilon \, .
\]
We have now shown that $D$ satisfies (i), (ii) and (iii) of Proposition \ref{wu-tr0}, whence $A$ has tracial rank zero. This completes the proof of Theorem \ref{lfdrtr0}.

\end{document}